\documentclass[fleqn,trsc,nonblindrev]{informs3noheader}

\OneAndAHalfSpacedXI 
\usepackage{endnotes}
\let\footnote=\endnote

\usepackage{graphicx}
\usepackage{multirow}
\usepackage{hhline}
\usepackage[english]{babel}
\usepackage[margin=1 in]{geometry}
\usepackage{verbatim}
\usepackage{multirow}
\usepackage{graphicx}
\usepackage[dvipsnames]{xcolor}
\usepackage{pgfplots}
\usepgfplotslibrary{colorbrewer}
\pgfplotsset{compat = 1.18, 
			 cycle list/Set1-8} 
\usetikzlibrary{pgfplots.statistics, pgfplots.colorbrewer} 
\usepackage{pgfplotstable}
\usepackage{filecontents}
\usepackage{textpos}

\usepackage{pifont}
\usepackage{threeparttable}
\usepackage{epstopdf}
\usepackage{booktabs}
\usepackage{makecell}
\usepackage[short, nocomma]{optidef}
\usepackage{natbib}
 \bibpunct[, ]{(}{)}{,}{a}{}{,}%
 \def\bibfont{\small}%
\usepackage{bbding}
\usepackage{csquotes}
\usepackage[colorlinks=true,linkcolor=blue,citecolor=blue,urlcolor=blue,]{hyperref}

\usepackage[small, margin=1cm]{caption}
\usepackage{subcaption}

\usepackage{appendix}

\usepackage{color}
\definecolor{strcolor}{rgb}{0.6, 0.2, 0.6}
\definecolor{commentcolor}{rgb}{0.3125, 0.5, 0.3125}
\definecolor{keycol}{rgb}{0, 0, 1}

\usepackage{bbm}

\usepackage{listings}
\usepackage{url}

\newcommand {\bea}{\begin{eqnarray}}
	\newcommand {\eea}{\end{eqnarray}}

\def\blot{\quad \mbox{$\vcenter{ \vbox{ \hrule height.4pt
				\hbox{\vrule width.4pt height.9ex \kern.9ex \vrule width.4pt}
				\hrule height.4pt}}$}}

\usepackage{natbib}
\bibpunct[, ]{(}{)}{,}{a}{}{,}%
\def\bibfont{\fontsize{8}{9.5}\selectfont}%

\TheoremsNumberedThrough    
\ECRepeatTheorems

\EquationsNumberedThrough   
\gdef\AQ#1{}
\gdef\CQ#1{}
\begin{document}
\def\COPYRIGHTHOLDER{INFORMS}%
\def\COPYRIGHTYEAR{2017}%
\def\DOI{\fontsize{7.5}{9.5}\selectfont\sf\bfseries\noindent https://doi.org/10.1287/opre.2017.1714\CQ{Word count = 9740}}

	\RUNAUTHOR{Lu et~al.} %

\RUNTITLE{Line planning under crowding: A cut-and-column generation approach}

\TITLE{Line planning under crowding: A cut-and-column generation approach}

	\ARTICLEAUTHORS{
		\AUTHOR{Yahan Lu,\textsuperscript{a, b} Rolf N. van Lieshout,\textsuperscript{c} Layla Martin,\textsuperscript{c, d} Lixing Yang\textsuperscript{a}} 
\AFF{$^{a}$School of Systems Science, Beijing Jiaotong University, Beijing, 100044, China; \\ $^{b}$Department of Transport \& Planning, Delft University of Technology, Netherlands; \\$^{c}$Department of Operations, Planning, Accounting, and Control, School of Industrial Engineering, Eindhoven University of Technology, Eindhoven, Netherlands; \\ ${^d}$ Eindhoven AI Systems Institute, Eindhoven University of Technology, Eindhoven, Netherlands}
}

\ABSTRACT{\textit{\textbf{Problem definition:}} To mitigate excessive crowding in public transit networks, network expansion is often not feasible due to financial and time constraints. Instead, operators are required to make use of existing infrastructure more efficiently. 
In this regard, this paper considers the problem of determining lines and frequencies in a public transit system, factoring in the impact of crowding. \textit{\textbf{Methodology:}} We introduce a novel formulation to address the line planning problem under crowding and propose a mixed-integer second-order cone programming (MI-SOCP) reformulation. Three variants of the cut-and-column generation algorithm with tailored acceleration techniques find near-system-optimal solutions by dynamically generating passenger routes and adding linear cutting planes to deal with the non-linearity introduced by the crowding terms. We find integral solutions using a diving heuristic. In practice, passengers may deviate from system-optimal routes. We, thus, evaluate line plans by computing a user-equilibrium routing based on Wardrop's first principle. \textit{\textbf{Results and implications:}} We experimentally evaluate the performance of the proposed approaches on both an artificial network and the Beijing metro network. The results demonstrate that our algorithm effectively scales to large-scale instances involving hundreds of stations and candidate lines, and nearly 57,000 origin-destination pairs. We find that considering crowding while developing line plans can significantly reduce crowding, at only a minor expense to the travel time passengers experience. This holds both for system-optimal passenger routing and user-optimal passenger routing, which only differ slightly. 
} 

	
\maketitle
\section{Introduction}
\label{sec:introduction}
Public transit passengers not only value fast and direct connections but also comfortable travel. Particularly, passengers prefer not to travel along highly crowded trajectories. However, due to the quickly growing population in metropolitan areas -- the UN expects 70\% of the world's population of then 9.7 billion to live in metropolitan areas by 2050 \citep{un75} -- public transit becomes increasingly crowded. 
High set-up costs and times make it impossible to adapt to this growth simply by extending and densifying the network. Public transit operators can, however, work towards utilizing the given rail network better by optimizing their planning. Especially changes to the \textit{line concept}, the set of lines that are operated and their frequencies, can potentially alleviate crowding, as these directly influence the routes taken by passengers. Therefore, in this paper, we focus on incorporating crowding considerations within the line planning decision, i.e., within a given network, deciding which lines to operate at which frequencies while considering passengers' travel time and crowding costs.

In real-world public transit networks such as the Beijing metro, with hundreds of stations and direct rails, determining the line concept is a challenging problem since passengers often transfer between different lines. This requires explicitly modeling transfers between lines and all possible routes that a passenger could take. Considering crowding in line planning further exacerbates complexity since it introduces interactions between passengers: The more passengers use the same train, the more severe the passengers perceive crowding. This poses a significant challenge in modeling crowding since the \textit{perceived travel time}, which is defined as the sum of travel time and crowding penalty, is no longer linear. Additionally, passengers might not travel in a centrally-optimal fashion through the network of public transit lines but rather find their own route.

To tackle these challenges, we develop a novel line planning model based on the change-and-go network introduced in \cite{Schobel2005} to trace the routes passengers take through the network. We model the impact of crowding on a passenger's perceived travel time as a linear function of the utilization of a train, resulting in a crowding term that is quadratic in the number of passengers, and inversely proportional to the frequency of that particular train. Through a convex reformulation based on second-order cone constraints, we obtain a computationally tractable mathematical model. We further observe that incorporating crowding effects into line planning generates several routes with similar perceived travel times. This results in small gaps between the system-optimal and user-optimal line plans, indicating that integrating user-optimal routing in line planning offers limited benefits. Therefore, we omit passengers' selfish route choices during the line planning optimization process but evaluate their impact \textit{ex-post}.

To solve the proposed line planning model at scale, we initially relax the cone constraints and start with a small subset of passengers' routes. Next, linear cuts (rows) are added if the cone constraints are violated, while routes (columns) are added if we find any negative reduced-cost routes through a pricing problem. Regarding the dynamic addition order of rows and columns, we develop two additional procedures: one that adds columns before rows and another that always adds columns and rows in each iteration. We obtain integral solutions using a tailored diving heuristic, which iteratively fixes lines with the highest frequency. Additionally, we develop acceleration techniques: We employ Farkas duals to deal with the infeasibility of the restricted master problem and remove non-tight cutting planes and columns with positive reduced costs where necessary. 

We test our approach using real-life instances from the entire Beijing metro network, one of the largest and most crowded systems in the world, consisting of 534 candidate lines and 56,916 origin-destination (OD) pairs. To create smaller instances, we extract the central subnetwork, which includes 168 candidate lines and 23,374 OD pairs. Additionally, we evaluate our approach on the most common benchmark in the line planning literature, the  $5 \times 5$-grid network \citep[e.g.,][]{SchobelAnita2017, Alexander2019, Schiewe2020}, where we generate a line pool with 128 lines. For the instances based on the grid network, the best procedure of our algorithm achieves an average optimality gap of less than 1.8\% compared to the optimal solutions obtained by GUROBI, while reducing the average computational time by more than 98.6\%. When solving the instance for the Beijing metro entire network without crowding effects, we are able to obtain the solution within 2.2 hours using the proposed algorithm. We solve instances on both the Beijing subnetwork and the entire network under crowding, obtaining higher-quality solutions that reduce passengers' perceived travel times compared to the existing line plan used in real-world operations. 

Extensive computational studies provide valuable insights into the benefits of incorporating crowding effects in line planning. In particular, the results show that incorporating crowding effects in line planning considerably reduces crowding, even when the crowding factor is small. Moreover, by incorporating crowding effects into line planning, passengers’ perceived travel time is reduced with only minor increases in travel time. Lastly, we show that the system-optimal routing only deviates slightly from the user-equilibrium, implying that the impact of ignoring passengers' selfish route choices during the optimization is limited.

To summarize, our contributions are four-fold: 

\begin{enumerate}
\item We introduce the line planning problem under crowding, which models crowding using a quadratic penalty term.  

\item Through a convex reformulation, we obtain a computationally tractable mixed-integer second-order cone program (MI-SOCP), which we solve using a cut-and-column generation framework.

\item We demonstrate that our algorithm obtains near-optimal solutions in a short computation time compared to GUROBI and is able to solve real-life instances.

\item From a practical standpoint, extensive numerical studies provide various managerial insights into line planning problems.

\end{enumerate}

The remainder of this paper is organized as follows. In Section \ref{sec:literature}, we review the literature related to this study. Section \ref{sec:problemDescription} formally introduces the line planning problem under crowding and presents the mathematical formulation. Subsequently, we present our solution approach, discussing the reformulation, the cut-generation method, and the cut-and-column generation algorithm in Section \ref{sec:algorithm}. Section \ref{sec:evaluateUE} proposes a model to evaluate the user equilibrium given the line plan. In Section \ref{sec:caseStudy}, we discuss the numerical results, including those based on the grid network and the Beijing metro system. Lastly, Section \ref{sec:conclusion} concludes this paper. 
\section{Literature review}
\label{sec:literature}

This paper merges two streams of literature on public transit systems: line planning that considers passenger routing and the effects of crowding on the public transit system. Section~\ref{sec:literLine} reviews line planning and describes how passenger routing can be integrated. We refer to \cite{Schobel2012} and \cite{schmidt2024planning} for a broader overview of the line planning literature. 
Section~\ref{sec:literCrowding} describes how existing literature models the effects of crowding within the planning phases of public transit.

\subsection{Line planning with passenger routing}
\label{sec:literLine}
The well-studied line planning problem usually selects several lines from a given line pool \citep{Gattermann2017}. Early contributions in line planning models exclusively focus on physical infrastructure and operational aspects \citep[e.g.,][]{Goossens2004,Goossens2006}. The goal is to find a line concept that serves all customers and minimizes operating costs for public transit companies.

More recent models have shifted towards a stronger focus on the passengers' preferences and objectives. Most of the related literature focuses on line planning with system-optimal routes, in parts integrated with other decisions such as network design \citep{Canca2019}, timetabling \citep{Lieshout2020}, or the integration of timetabling and vehicle scheduling \citep{Schobel2017}. In these problems, the operator chooses lines that optimize their objective (e.g., profit, costs, service level) by assigning passengers to routes. This stream of literature can be divided into two substreams. The first substream maximizes the number of passengers traveling without transfers. Those who must make transfers to reach their destinations are often overlooked. For example, the first model in \cite{Bertsimas2021} takes passengers with direct routes into account and maximizes the total demand that can be served by the line plan, where capacity limitations are disregarded. In this model, passenger preferences are limited to either choosing a direct route or forgoing being served.

The second substream minimizes the total travel time of passengers, taking routes with transfers into account. Operating costs usually translate to budget constraints, while capacity limitations are either discarded or formulated as hard constraints.  
Tracking passengers' travel time on each available route, \cite{Schobel2005} propose the change-and-go network (CGN). This network representation allows for characterizing passenger transfer behavior precisely and formulating passengers' travel time as a function of both in-vehicle time and transfer penalties. For example, \cite{SchobelAnita2017} present a line planning model incorporating hard capacity constraints, wherein passenger assignment is integrated into cost-optimal line planning. This model ensures that the average length of the considered routes does not exceed a specified percentage increase compared to the shortest route. The authors test this model on the $5 \times 5$-grid network with up to 275 lines in the line pool. As an alternative to modeling based on the CGN, \cite{Yao2024} propose a two-phase solution approach for the railway line planning problem with passenger routing using a direct-connection network. In this model, the capacity of lines and trains is ignored and the requirement to serve all passengers is relaxed. The line planning model proposed in \cite{Bertsimas2021} uses yet another modeling alternative, which considers routes with at most one transfer. In this model, the travel time for routes with a transfer is constrained by an exogenous threshold. Routes are all pre-computed outside the optimization process. Moreover, this model discards capacity limitations and aims to maximize the number of demands served by the transit network.

Beyond the aforementioned approach of pre-computing passengers' routes outside of the optimization framework, \cite{Borndorfer2007} use a column generation algorithm to dynamically generate passenger routes in the physical network within the optimization process. They develop a mixed-integer model aimed at minimizing the sum of passengers' travel costs and operators' costs. However, this model excludes transfer penalties and overlooks capacity limitations on lines. 

More recently, the line planning literature moved away from assigning passengers to routes and rather explicitly models passenger route choice using discrete choice models and Nash equilibria \citep{Schmidt2012, Hartleb2023}. For example, \cite{Goerigk2017} propose a bi-level model based on the CGN for the line planning problem that incorporates user-optimal route choices. In this model, passenger routes are generated prior to optimization, and capacity limitations are explicitly formulated as hard constraints. The travel time of passengers is formulated as a function associated with in-vehicle travel time and transfer penalties. The authors develop a genetic algorithm to solve large-scale instances. \cite{Alexander2019} interpret line planning as a game in which passengers are players striving to minimize their individual objective functions, which comprise travel time, transfer penalties, and a share of the overall solution cost, while accounting for hard capacity constraints. Both methods face numerical difficulties due to the interaction terms' bi-level nature and non-convexity. 

While the line planning literature has moved towards more detailed modeling of passenger preferences by either minimizing or constraining the overall travel time or transfers, passengers are not only sensitive to physical aspects of their journey but also travel comfort, such as crowding. We thus integrate the effects of crowding on the perceived travel time of passengers in the line planning problem and model the impact of crowding as soft-capacity constraints. In addition, current algorithms addressing line planning problems exhibit limitations in scalability. Consequently, we develop a cut-and-column generation approach designed to solve real-life instances involving hundreds of stations, hundreds of candidate lines in the line pool, and tens of thousands of OD pairs.

\subsection{Crowding effects in public transit planning}
\label{sec:literCrowding}
More passengers on the same train reduces travel comfort. This has not been considered in the context of line planning and is usually modeled using hard capacity constraints in other public transit decisions, such as timetabling \citep[e.g.,][]{Cacchiani2020, Yin2021, Lu2022, Lu2023}, vehicle scheduling \citep[e.g.,][]{Xia2024integrated, Chai2024}, and charging \citep[e.g.,][]{Marelot2024}. 

The literature that precisely models the effects of crowding on passengers' travel discomfort in public transit planning process is quite limited. \cite{Jiang2022} formulate a mixed integer nonlinear programming (MINLP) model for the timetabling problem in a transit network offering both schedule- and frequency-based services. This model incorporates different in-vehicle travel times for passengers, whether seated or standing, using a piecewise linear function. \cite{Luan2022} present a MINLP formulation to reschedule trains and reroute passengers, incorporating crowding effects into system-optimal passenger movements. They assume that crowding effects do not come into play if every passenger has a seat but increase stepwise otherwise. These effects are modeled as a piecewise constant function within the formulation.  

The freight transport literature considers crowding more frequently. More vehicles on the same road increase the travel time \citep[e.g.,][]{Daganzo1997, Franceschetti2018}, drivers reroute accordingly, and thus road network planners include crowding. Similarly, network design and location decisions in centrally planned logistics systems benefit from considering crowding if demand is stochastic \citep[e.g.~][]{Vedat2023,Burak2023}. For example, \cite{Basciftci2023} consider the crowding effects on hub location decisions, factoring in strict capacity limitations. The authors employ the Kleinrock function to model each hub as an M/M/1 queue under steady-state conditions, thereby formulating the crowding effects as a nonlinear cost function.

We add to this literature stream by considering crowding in a decision-making problem that previously ignored these effects. We propose a novel crowding effect function and develop an efficient cutting-plane method to address the computational challenges arising from the quadratic crowding terms. Additionally, we design three procedures within the cut-and-column generation approach to efficiently handle the quadratic crowding terms when dealing with an enormous number of columns simultaneously.
\section{Problem statement}
\label{sec:problemDescription}

In this paper, we consider the \textit{line planning problem under crowding} (LPP-C). 
We first present the underlying public transit network (PTN) and then introduce the corresponding change-and-go network (CGN) in line with \cite{Schobel2005}. We subsequently explain the LPP-C in greater detail. Lastly, we formulate the problem as a mathematical model with a nonlinear objective function.

\subsection*{Network}
A PTN is a graph $G = (S, E)$ with a set of vertices $S$ representing stations and a set of edges $E$ representing tracks between two stations. A line is a simple path in $G$, and all potential lines are given by the line pool $\mathcal{L}$. This line pool is a subset of all simple paths. Each station should be served by at least one line; if multiple lines pass a station $i\in S$, this station is a possible transfer location. 

To track passengers changing between lines and account for the inconvenience associated with transfers, we cast the PTN into a CGN. A CGN is a graph $\mathcal{G} = (\mathcal{N}, \mathcal{A})$, where $\mathcal{N}$ refers to the set of nodes in the CGN. The node set $\mathcal{N}$ comprises the set of station nodes $S$ and the set of travel nodes $N^\text{t}$. We generate a node in the CGN for each station-line $(s,l)$ tuple such that $s\in l$ allows us to track passengers transferring between lines explicitly. An example of a four-station PTN and its corresponding CGN with four station nodes and six transfer nodes is given in Figure \ref{fig_examplePTN&CGN}. The set $\mathcal{A}$ contains travel arcs $\mathcal{A}^\text{t}$ and transfer arcs $\mathcal{A}^\text{c}$. Each travel arc $a \in \mathcal{A}^\text{t}$ incurs a perceived cost comprising the travel time and a crowding penalty. The latter negatively correlates with frequency and positively depends on the utilization. In contrast, transfer arcs only incur a fixed cost $c_a$ for making a transfer.

\input{Tables_and_figures/fig_examplePTN_CGN}

\subsection*{The line planning problem under crowding}
As in most line planning models, we assume the line pool is defined a priori \citep{Lieshout2020}. From the line pool $\mathcal{L}$ including a set of candidate lines connecting stations in a public transit system, the operator determines which lines to operate at what frequency. The operator aims to minimize passengers' perceived travel times by operating lines subject to a budget constraint $B$ such that all demand is served.

Once the lines have been established, passengers join the system and travel from origin to destination. We thus aggregate interchangeable passengers with the same origin and destination. We denote the passenger demand by the set $\mathcal{P}\subseteq S \times S$, which includes the passengers' origin and destination, and use $n^p$ for the positive passenger demand of $p\in \mathcal{P}$. Every passenger is assigned to a route, which is a directed path in the CGN that starts from the station node corresponding to the passengers' origin, passes through various nodes and arcs, and arrives at the station node corresponding to the destination. Along this route, passengers can transfer between different lines. We incorporate crowding effects into line planning, minimizing a perceived travel time function of passengers comprising travel time and crowding penalty. 

\subsection*{Mathematical formulation}
To formulate the LPP-C, we first introduce a mathematical model with a nonlinear objective function. We introduce the parameter $\mathcal{L}(e)$ to represent the set of lines passing through edge $e \in E$ in the PTN. Let $\mathcal{R}$, $\mathcal{R}_a$, $\mathcal{R}_p$, $\mathcal{R}_a^p$ denote the set of possible routes for all OD pairs, the set of routes passing arc $a \in \mathcal{A}$, the set of possible routes for OD pair $p\in \mathcal{P}$, and the set of routes for OD pair $p$ which pass arc $a$. Let $d_{\ell}$ and $e_{\ell}$ denote the operational costs and set-up costs of line $\ell$, respectively. Let $\ell(a)$ represent the line $\ell$ to which arc $a$ belongs. We define the parameter $\gamma_a>0$ to represent the \textit{crowding factor}, which depends on the capacity of the vehicle that operates on $\ell(a)$. The budget is denoted as $B$, and the minimum and maximum frequency of line $\ell$ are denoted as $f^{min}_{\ell}$ and $f^{max}_{\ell}$, respectively. Four types of decision variables represent the line plan and the passenger flows:

(i) $y_{\ell}$: Continuous variables representing the frequency of line $\ell \in \mathcal{L}$,

(ii) $w_{\ell}$: Binary variables indicating whether line $\ell \in \mathcal{L}$ is operated in the solution or not,

(iii) $x_a$: Continuous decision variables representing the total number of passengers using arc $a \in \mathcal{A}$,

(iv) $z_r$: Continuous decision variables indicating the number of passengers that choose route $r \in \mathcal{R}$. 

Then, the LPP-C can be formulated as follows:
\begin{mini!}[0]
{}
{\sum_{a\in \mathcal{A}}c_{a}x_{a} + \sum_{a\in \mathcal{A}: y_{\ell(a)>0}}\frac{\gamma_a x_{a}^2}{y_{\ell(a)}} \label{eq:objective}}
{\label{formulation:NL}} 
{} 
\addConstraint
{\sum_{\ell \in \mathcal{L}}\left(d_{\ell} y_{\ell} + e_{\ell}w_{\ell}\right)}
{\leq B \quad \label{eq:budgetConstraint}}
{}
\addConstraint
{\sum\limits_{\ell \in \mathcal{L}(e)} w_{\ell}}
{\geq 1 \label{eq:connectedConstraint}}
{\ \forall e \in E,}
\addConstraint
{f^{min}_{\ell}w_{\ell}}
{\leq y_{\ell} \leq f^{max}_{\ell} w_{\ell} \quad \label{eq:lineNumberConstraint1}}
{\ \forall \ell \in \mathcal{L},}
\addConstraint
{x_{a}}
{=\sum_{r\in\mathcal{R}_{a}}z_r \quad \label{eq:SOflow_arcConstraint}}
{\ \forall a \in \mathcal{A},}
\addConstraint
{\sum_{r\in \mathcal{R}_a^p}z_r}
{\leq n^p w_{\ell(a)} \quad \label{eq:SOflow_arc_lineConstraint} }
{\ \forall  a\in\mathcal{A}, p \in \mathcal{P},}
\addConstraint
{\sum_{r\in \mathcal{R}_p}z_r}
{=n^p \quad \label{eq:SOProbabilityConstraint}}
{\ \forall p \in \mathcal{P},}
\addConstraint
{w_{\ell}}
{\in \{0, 1\} \label{eq:lineBinary_domainConstraintNew}}
{\ \forall \ell \in \mathcal{L},}
\addConstraint
{x_a}
{ \geq 0  \quad \label{eq:SOflow_arc_od_domainConstraint}}
{\ \forall a \in \mathcal{A},}
\addConstraint
{y_{\ell}}
{\geq 0 \quad \label{eq:line_domainConstraint}}
{\ \forall \ell \in \mathcal{L},}
\addConstraint
{z_r}
{\geq 0 \quad \label{eq:SOroute_od_domainConstraint}}
{\ \forall r \in \mathcal{R}.}
\end{mini!}
The objective function \eqref{eq:objective} minimizes the total perceived travel times of passengers. Constraint \eqref{eq:budgetConstraint} ensures that the operational and set-up costs of the operated lines remain within the budget. Constraints~\eqref{eq:connectedConstraint} ensure that every edge in the PTN has at least one line passing through it. Constraints~\eqref{eq:lineNumberConstraint1} link variables $y_{\ell}$ and $w_{\ell}$, enforcing that operated lines have a frequency between $f^{\min}_{\ell}$ and $f^{\max}_{\ell}$ and 0 otherwise. Constraints~\eqref{eq:SOflow_arcConstraint} link the arc flows to the route variables, making sure that the flow through an arc equals the flow over all routes passing through this arc. Constraints \eqref{eq:SOflow_arc_lineConstraint} restrict passengers to use arcs that belong to operated lines $\ell$, i.e., for all $\ell$ with $w_{\ell} = 1$. Constraints~\eqref{eq:SOProbabilityConstraint} ensure that all passengers are routed. Constraints \eqref{eq:lineBinary_domainConstraintNew} – \eqref{eq:SOroute_od_domainConstraint} define the domains of the decision variables.

\section{Solution methodologies}
\label{sec:algorithm}

Formulation \eqref{formulation:NL} cannot be used to solve realistic instances using standard commercial solvers because (i) the set over which we sum in the crowding penalty term in the objective function depends on the values of the decision variables and (ii) the number of routes for every OD pair may be enormous. To deal with (i), Section~\ref{subsec:reform} 
reformulates the problem as a mixed-integer second-order cone program (MI-SOCP). Section~\ref{sec:cuts} describes how this MI-SOCP can be solved using a cutting-plane approach. To handle (ii), Section~\ref{sec:columnGeneration} shows how routes can be generated dynamically using column generation. Section~\ref{sec:overviewAlgorithm} introduces three variants of the proposed cut-and-column generation approach, which integrates the two algorithms described above. Lastly, Section~\ref{sec:diving} presents the diving heuristic to obtain integral solutions.

\subsection{Reformulation as a mixed-integer second-order cone program}
\label{subsec:reform}

To reformulate the problem as an MI-SOCP, we introduce auxiliary decision variables $\Theta_a$ for all $a \in \mathcal{A}$ that serve as the crowding penalty incurred on arc $a$ in the objective function. To ensure that these variables take the right value, we add a set of rotated second-order cone constraints. In the reformulation, we no longer need constraints \eqref{eq:SOflow_arc_lineConstraint} to couple decision variables $\mathbf{z}$ and $\mathbf{w}$. This results in the following formulation: 
\begin{mini!}[0]
{}
{\sum_{a \in \mathcal{A}} c_{a}x_{a}+ \sum_{a \in \mathcal{A}} \Theta_a \label{eq:objTightening}} 
{\label{formulation:SOconvexStrong}} 
{}    
\addConstraint
{\sum_{\ell \in \mathcal{L}}\left(d_{\ell} y_{\ell} + e_{\ell}w_{\ell}\right)}
{ \leq B  \quad \label{eq:budgetConstraintNew}}
{}
\addConstraint
{\Theta_a y_{\ell(a)}}
{\geq \gamma_a x_a^2 \quad \label{eq:tightLinedemand}}
{\ \forall a \in \mathcal{A},}
\addConstraint
{\sum\limits_{\ell \in \mathcal{L}(e)} w_{\ell}}
{\geq 1 \label{eq:connectedConstraint1}}
{\ \forall e \in E,}
\addConstraint
{f^{min}_{\ell}w_{\ell}}
{\leq y_{\ell} \leq f^{max}_{\ell} w_{\ell} \quad \label{eq:lineNumberConstraint1New}}
{\ \forall \ell \in \mathcal{L},}
\addConstraint
{x_{a}}
{=\sum_{r\in\mathcal{R}_{a}}z_r \quad\label{eq:SOflow_arcConstraint2}}
{\ \forall a \in \mathcal{A},}
\addConstraint
{\sum_{r\in \mathcal{R}_p}z_r}
{=n_p \label{eq:SOProbabilityConstraint2}}
{\ \forall p \in \mathcal{P},}
\addConstraint
{w_{\ell}}
{\in \{0, 1\} \label{eq:lineBinary_domainConstraintNew2}}
{\ \forall \ell \in \mathcal{L},}
\addConstraint
{x_a}
{ \geq 0  \quad \label{eq:SOflow_arc_od_domainConstraint2}}
{\ \forall a \in \mathcal{A},}
\addConstraint
{y_{\ell}}
{\geq 0 \label{eq:line_domainConstraint2}}
{\ \forall \ell \in \mathcal{L},}
\addConstraint
{z_r}
{\geq 0 \label{eq:SOroute_od_domainConstraint2}}
{\ \forall r \in \mathcal{R},}
\addConstraint
{\Theta_a}
{\geq 0 \label{eq:theta_domainConstraint}}
{\ \forall a \in \mathcal{A}.}\end{mini!}

Constraints \eqref{eq:tightLinedemand} are the most important constraints in Formulation~\eqref{formulation:SOconvexStrong}. If $y_{\ell(a)} = 0$, these constraints force $x_a = 0$ and if $y_{\ell(a)}>0$, they ensure that $\Theta_{a}\geq \gamma_a x_a^2/{y_{\ell(a)}}$. These constraints are nonlinear, but, together with the conditions that $y_{\ell(a)}\geq 0$ and $\Theta_a\geq 0$ for all $a$, describe a rotated second-order cone, which is convex. The remaining constraints follow from the original formulation \eqref{formulation:NL}.

\subsection{Cutting planes for SOCP constraints}
\label{sec:cuts}

Mixed-integer second-order cone programs can be solved directly by most contemporary commercial solvers. As a faster alternative, the rotated second-order cone constraints can be handled using cutting planes that also integrate well with column generation.
In this method, the cone constraints are relaxed, yielding a linear formulation, and tangent cuts are added on-the-fly when violations of the constraints are detected, as formalized in the following proposition:

\begin{proposition}[Cut Generation]
\label{lemma:cut}
Given a point ($\hat{\Theta},\hat{x},\hat{y}$) where $\hat{y}_a > 0$ for some $a\in \mathcal{A}$, the following cut is valid:
\begin{align}
   \frac{\Theta_{a}}{\gamma_a}  \geq \frac{2\hat{x}_a}{\hat{y}_{\ell(a)}}x_a-\frac{\hat{x}_a^2}{\hat{y}^2_{\ell(a)}}y_{\ell(a)}.\label{cut1}
\end{align}
\end{proposition}

\proof{Proof.}
We will show that this cut does not eliminate any feasible solutions of Formulation~(\ref{formulation:SOconvexStrong}). We omit indices for brevity. Let $\left(\Theta,x,y\right)$ denote a feasible solution. If $y=0$, it must be that $\Theta=x=0$, and (\ref{cut1}) is satisfied. If $y>0$, constraint \eqref{eq:tightLinedemand} becomes $\Theta/\gamma \geq  f(x,y)$, where $f(x,y) =  x^2/y$. Since $f(x, y)$ is convex, given the point $(\hat{x},\hat{y})$, it holds that 
$$f(x, y) \geq f(\hat{x},\hat{y})+ f_x(\hat{x},\hat{y}) (x-\hat{x})+ f_y(\hat{x},\hat{y}) (y-\hat{y}),$$
where $f_{\cdot}$ denotes the subgradient. Therefore, we have
\begin{align*}
\frac{\Theta}{\gamma}\geq f(x, y) 
 & \geq \frac{ \hat{x}^2}{\hat{y}} +  \frac{2  \hat{x}}{\hat{y}} (x-\hat{x}) -\frac{ \hat{x}^2}{\hat{y}^2}  (y-\hat{y}) \\
& \geq \frac{2  \hat{x}}{\hat{y}}x-\frac{ \hat{x}^2}{\hat{y}^2}y, 
\end{align*}
which corresponds to \eqref{cut1}.
\Halmos
\endproof

Our algorithm checks for constraint violation and adds these linear cuts when necessary. In case $y_l=0$ for some violated constraint, we replace $y_l$ by $y_l+\varepsilon$ for some small $\varepsilon$ before computing the cut. 

\subsection{Column generation}
\label{sec:columnGeneration}

The other challenge of solving Formulation~\eqref{formulation:SOconvexStrong} lies in the enormous set of routes $\mathcal{R}_{p}$ for each OD pair $p \in \mathcal{P}$ when the network is large, which results in a huge number of variables. However, the challenge can be handled efficiently using the column generation algorithm. 

\subsubsection{Iterative process of the column generation algorithm} Our column generation algorithm proceeds iteratively as follows. First, we begin with a restricted set of routes $\overline{\mathcal{R}} \subset \mathcal{R}$ and the corresponding subset of variables. We solve the \textit{restricted master problem} (RMP), which is Formulation \eqref{formulation:SOconvexStrong} on $\overline{\mathcal{R}}$ where we additionally relax the integer condition on $\mathbf{w}$ and remove constraints \eqref{eq:tightLinedemand} to obtain an optimal primal and dual solution to the RMP. Then, to see if the solution is also optimal with respect to the full set of routes $\mathcal{R}$, we solve the \textit{pricing problem} (PP) by searching for the route $r\in \mathcal{R}$ with the smallest reduced cost given the current dual variables. For any given route $r$ for OD pair $p$, the reduced cost $\psi^p_r$ can be expressed as $\psi^p_r =-\lambda^p+\sum\limits_{a\in \mathcal{A}_r}\vartheta_a$,
where $\vartheta_a$ is the dual of constraints \eqref{eq:SOflow_arcConstraint2}, and $\lambda^p$ is the dual of constraints \eqref{eq:SOProbabilityConstraint2}. For each $p\in \mathcal{P}$, we can find the route $r^*$ that minimizes the reduced cost by solving a shortest-path problem on a directed acyclic graph. If $\psi^p_{r^*}<0$, adding $r^*$ to $\overline{\mathcal{R}}$ can improve the current solution. If $\psi^p_{r^*}\geq 0$ for all $p\in P$, the current solution is optimal. 

\subsubsection{Initialization and dealing with infeasibility}
To initialize the column generation algorithm, we include the shortest route in the CGN for each OD pair in $\overline{\mathcal{R}}$. Given the restricted number of routes in $\overline{\mathcal{R}}$, the RMP may become infeasible. This issue would also arise in the process of finding integral solutions due to the possible absence of routes for certain OD pairs in the line plan with all integral frequencies. The infeasibility arises from constraints~\eqref{eq:SOProbabilityConstraint2} not being satisfied and the associated dual problem being unbounded. During column generation, utilizing Farkas pricing enables the restoration of the feasibility of the RMP \citep{Kowalczyk2018}. To be more precise, we retrieve the Farkas duals of the constraints \eqref{eq:SOflow_arcConstraint2} and \eqref{eq:SOProbabilityConstraint2}, use them to solve the pricing problem, add routes to the RMP, and solve the RMP again. We repeat these procedures until we have shown that the RMP is feasible.

\subsection{An overview of the algorithm}
\label{sec:overviewAlgorithm}

We investigate three variants of the cut-and-column generation approach: 

(i) First-Cut-Then-Price (FCTP). After solving the RMP, this method first identifies and adds violated cuts, and only calls the pricing problem if no cuts are added. 

(ii) First-Price-Then-Cut (FPTC). Conversely, this method always calls the pricing problem after solving the RMP, and only adds cuts if there are no negative reduced cost columns. 

(iii) Always-Price-And-Cut (APAC). The final method always calls both the pricing and the cutting procedure after solving the RMP. 

To limit the size of the models, we apply both cut and column management. In every proposed variant, each time the cut generation method is called, we identify the tight cuts, save them for the next iteration, and eliminate the slack ones so as to further speed up the computation. Additionally, after each pricing procedure, we compute the average number of routes with respect to all OD pairs. If this measure exceeds a predetermined threshold (denoted as $\Delta_1$), we remove routes whose reduced costs exceed a positive threshold value $\Delta_2$. 

Termination options for this heuristic involve solving the RMP to optimality or stopping the column generation process before the optimal solution of the RMP is reached to avoid the tail-off effect. Here, we terminate the column generation process if the relative improvement in the objective value of the RMP over the last $I$ iterations falls below $\epsilon$. These parameters dynamically adapt during the solution process. 

\subsection{Diving heuristic}
\label{sec:diving}

To obtain integral solutions, we utilize the diving heuristic to perform a heuristic search in the branch-and-bound tree. The diving heuristic iteratively executes a column generation phase and a fixing phase until it yields an integral solution. Specifically, after a column generation phase, we round up the fractional line variable $w_{\ell}$ whose corresponding frequency variable $y_{\ell}$ is the largest. Thereafter, we check whether the remaining budget is sufficient to operate another line at its minimum frequency. If the budget is sufficient, we continue the search by solving the RMP again. If not, the diving heuristic terminates. 

Since we use an early stopping criterion to solve the RMP, it is possible that there are still negative reduced cost routes that can improve the objective, whilst keeping the line variables fixed. Therefore, after terminating the diving heuristic, we reoptimize the RMP and resolve the pricing problem until no negative reduced cost routes can be found.

\section{Evaluating the user equilibrium}
\label{sec:evaluateUE}
Given a line plan, we are interested in evaluating the user equilibrium considering the crowding effects. The user equilibrium principle of Wardrop states that each passenger strives to minimize their individual perceived travel times and will collectively reach an equilibrium state where no passenger has an incentive to change their chosen route. Given a solution $\mathbf{y}^*$, let $\mathcal{A}^{*}$ denote the set of all arcs $a$ such that $y_{l(a)}^*>0$, and define the function $d_a \left(x_a \right):= c_a+ \frac{\gamma_a}{y_{l(a)}}x_a$ for the \textit{perceived travel cost} on each arc $a \in \mathcal{A}^{*}$. To formally define the equilibrium conditions, let the continuous variable $\tau^{p}$ represent the travel costs in the equilibrium state of each OD pair $p \in \mathcal{P}$ and let the continuous variable $h_{r}$ represent the travel cost of each route $r$ in $\mathcal{R}$. A routing $\mathbf{z}$ is an equilibrium if $(\mathbf{x},\mathbf{\tau})$ exist such that the following conditions are satisfied \citep{Wardrop1952}:
\begin{align*}
\sum\limits_{r\in\mathcal{R}_{p} }z_{r} & = n_p & \forall p \in \mathcal{P}, \\
x_a & = \sum\limits_{r\in\mathcal{R}_{a} }z_{r}  & \forall a \in \mathcal{A}^{*}, \\
h_{r} &= \sum_{a\in r} d_a \left(x_a \right) & \forall  r\in\mathcal{R},\\
z_r > 0 &\Rightarrow h_{r} = \tau^{p} & \forall p \in \mathcal{P},  r \in \mathcal{R}_p, \\
z_r = 0 &\Rightarrow h_{r} \geq \tau^{p} &\forall p \in \mathcal{P}, r \in \mathcal{R}_p.
\end{align*}
These conditions are non-linear and, therefore, hard to solve directly. However, as shown in Proposition 18.11 in \cite{roughgarden2007}, we can find an equilibrium routing by solving the following optimization problem:
\begin{mini!}[2]
{}
{\sum_{a \in \mathcal{A}^{*}}\int_{0}^{x_a} d_a(s)ds }
{\label{formulation:classicUE}} 
{}
\addConstraint
{x_{a}}
{=\sum_{r\in\mathcal{R}_{a}}z_r }
{\ \forall a \in \mathcal{A}^{*},}
\addConstraint
{\sum_{r\in \mathcal{R}_p}z_r}
{=n^p }
{\ \forall p \in \mathcal{P},}
\addConstraint
{x_a}
{ \geq 0  }
{\ \forall a \in \mathcal{A}^{*},}
\addConstraint
{z_r}
{\geq 0 }
{\ \forall r \in \mathcal{R}.}
{}
{}
\end{mini!}
For any continuous and non-decreasing function $d_a(x_a)$, the objective function  is both continuously differentiable and convex. For our case, the integral evaluates to:
\begin{align}
    \sum\limits_{a \in \mathcal{A}^{*}} \int_{0}^{x_a} d_{a}(s) ds     & = \sum\limits_{a \in \mathcal{A}^{*}}\int_{0}^{x_a} \left(c_{a}+ \frac{\gamma_a}{y^*_{\ell(a)}} s \right) ds\nonumber \\
    & =\sum\limits_{a \in \mathcal{A}^{*}} \Bigg[c_a s + \frac{\gamma_a}{2y^*_{\ell(a)}} s^2 \Bigg]_{0}^{x_a}  \nonumber \\
    & =\sum\limits_{a\in \mathcal{A}^{*}}c_{a}x_{a}+\sum\limits_{a\in \mathcal{A}^{*}}\frac{\gamma_a}{2 y^*_{\ell(a)}} x_a^2.
\end{align} 
Therefore, for a given line plan, we can find its corresponding user equilibrium by solving a convex quadratic problem. Note that since the objective function and constraint functions are convex, the user equilibrium is unique. For cases where the route set $\mathcal{R}$ is very large, it can again be solved using column generation, using the same approach as outlined in Section~\ref{sec:columnGeneration}. 

\section{Numerical experiments}
\label{sec:caseStudy}

We now turn to computational experiments conducted on three networks using artificial and real-life operational data, respectively. The algorithm is implemented in Java, and all experiments are performed on a personal computer equipped with a 12th Gen Intel Core i7-12700H CPU (2.30 GHz), 64 GB of RAM, 16 cores, and GUROBI 9.5.1.

When analyzing the solutions to the instances, we aim to answer questions such as:

(i) How much does incorporating passenger routing help in reducing perceived travel times? When using GUROBI as a benchmark solution method, how many routes are appropriate to include in the optimization process?

(ii) How much do our three variants of the developed cut-and-column generation algorithm help in improving the solution efficiency in both artificial and real-life networks? What are the characteristics of the most applicable scenarios for each variant?

(iii) What are the characteristics of high-quality line plans when minimizing passengers' perceived travel times under crowding? What are the benefits of optimizing line plans compared to the existing ones used in current practical operations?

(iv) What are the benefits of incorporating crowding effects into line planning to enhance passengers' travel experience, particularly in reducing crowding?

(v) What are the characteristics of deviations between system-optimal and user-equilibrium routing when minimizing passengers' perceived travel times under crowding?

Before analyzing the results of our numerical experiments, we first provide details of the instances used in our study in Section \ref{sec:numericalDesign}. In Section \ref{sec:GridResults}, we address the first question regarding the effectiveness of including passenger routing and partially discuss the second question on the solution efficiency of the cut-and-column generation algorithm. Section \ref{sec:managerialInsights} focuses on the remaining three questions, exploring the benefits of optimizing line plans, the characteristics of high-quality solutions, the benefits of incorporating crowding effects, and deviations between system-optimal and user-equilibrium routing. In Section~\ref{sec:managerialInsights}, we also demonstrate that our methodologies generate more effective line plans relative to the existing line plan employed in real-world operations, leading to smaller passengers' perceived travel times within the network. Lastly, in Section \ref{sec:entireNetwork}, we demonstrate the scalability of our algorithm on the entire Beijing metro network, providing further insights into the second question on large-scale problem-solving capability.

\subsection{Experimental design}
\label{sec:numericalDesign}

To demonstrate the effectiveness of the proposed formulations and algorithms, we present the results of numerical experiments utilizing artificial instances with the grid network and real-world case studies based on the Beijing metro network. Table~\ref{tab_characteristics} summarizes the characteristics of the instances.

\begin{table}[h]
\caption{Characteristics of the instances.}
\label{tab_characteristics}
\centering
\begin{tabular}{cccrr}
\Xhline{1pt}
Instance & Network         & Candidate lines (\#) & OD pairs (\#) & Demand (\#) \\
\Xhline{0.6pt}
A        & $5 \times 5$-grid network   & 128  &  567  &   2,547        \\
B        & \begin{tabular}[]{@{}c@{}}Beijing metro central \\  sub-network with 8 existing lines \end{tabular}    &   168            &             23,374      &   446,649         \\
C        & Beijing metro entire network    &       534        &       56,916              &    776,401           \\
\Xhline{1pt}
\end{tabular}
\end{table}

First, we focus on a $5 \times 5$-grid network detailed in \cite{GridRef} and \cite{schiewe2024lintim} to analyze the benefits of incorporating passenger routing into the line planning problem and the effectiveness of our algorithm. As Figure~\ref{fig_gridNetwork} shows, this network comprises 25 nodes and 40 edges, and serves 2,547 passengers in 567 OD pairs. Travel times on each edge are 6 minutes and dwell times at stations are 1 minute. We construct a line pool consisting of 128 lines. The operational costs $d_\ell$ for line $\ell$ is equal to $8N_\ell$ where $N_\ell$ is the number of edges traversed by line $\ell$. Additionally, the set-up costs $e_{\ell}$ are positively proportional to 
$d_\ell$, and are given by $e_{\ell} = 1.5 \times d_\ell$. We obtain a baseline budget $\underline{B}$ by solving the cost-minimal line planning problem with strict capacities. In the experiments, we vary the actual budget $B$ to 100\% and 120\% of this value. The transfer penalty is set at 5 minutes. Based on a large number of pre-experiments with $B$ set at $\underline{B}$, we observe that the maximum frequency of operated lines in optimal solutions consistently remains no greater than 5, regardless of variations in the crowding factor. This observation supports the conclusion that higher frequency values are unnecessary in this context. Therefore, in all experiments based on the grid network, the upper bound for frequency (denoted as $f^{max}_{\ell}$ for all $\ell \in \mathcal{L}$) is set to 5 trains per hour.

\input{Tables_and_figures/fig_gridNetwork}

Thereafter, we perform experiments on the Beijing metro central sub-network to explore the benefits of the optimized line plans obtained by the algorithm over the line plan used in practical operations. Lastly, we turn to the entire Beijing metro network to showcase the capability of the proposed cut-and-column generation algorithm to solve large-scale problems. For these instances based on the Beijing metro sub-network and entire network, we set the maximum frequency to 40 trains per hour, as proposed by operational managers of the Beijing metro \citep{beijing2023subway}.

\begin{figure}[tbp]
    \centering
     \begin{subfigure}[b]{\textwidth}
         \centering
         \includegraphics[width=0.8\linewidth]{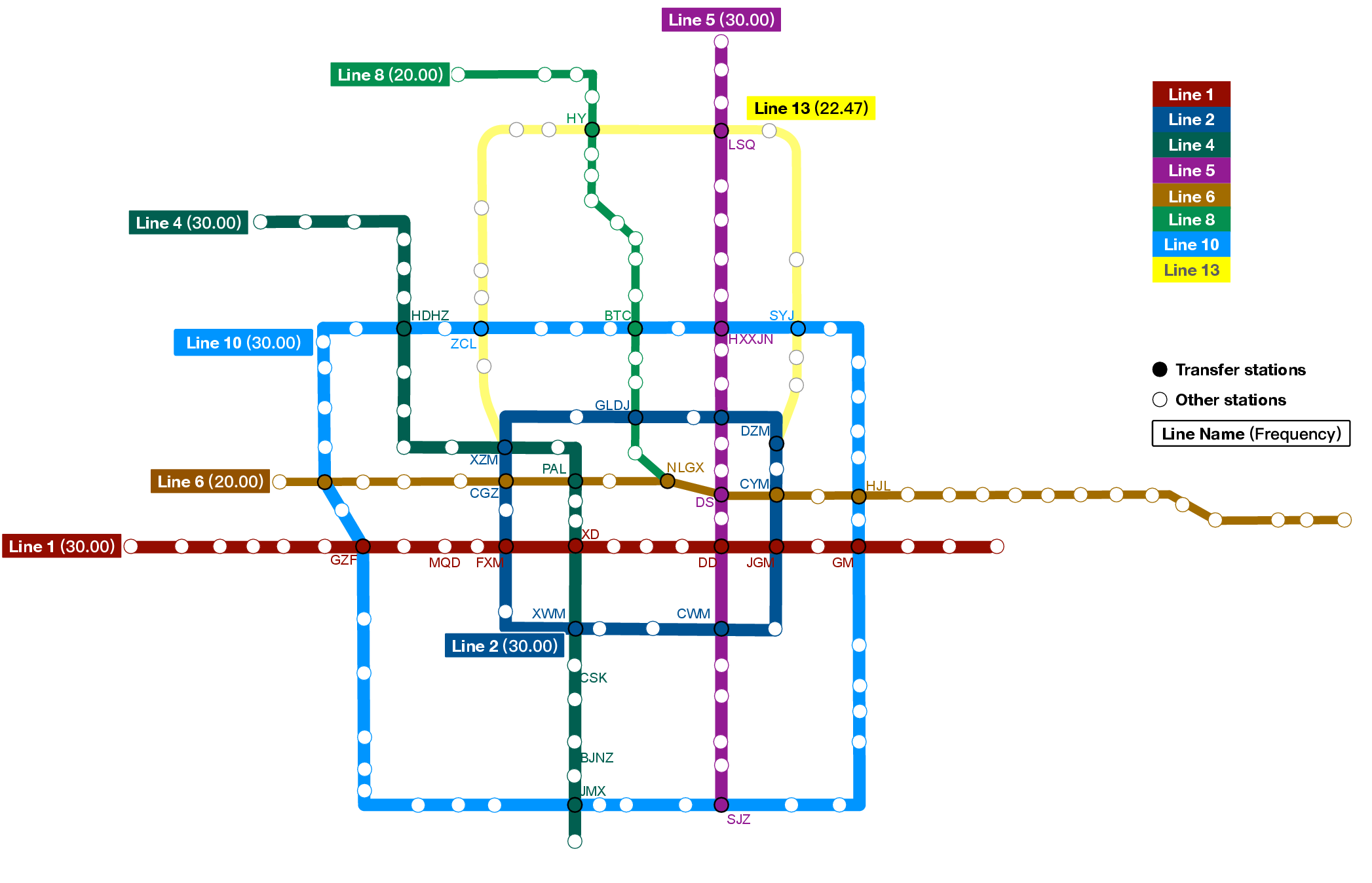}
         \caption{Existing line plan in the Beijing metro central sub-network}
     \end{subfigure}

    \vspace{1em}
    
     \begin{subfigure}[b]{\textwidth}
         \centering
         \includegraphics[width=\linewidth]{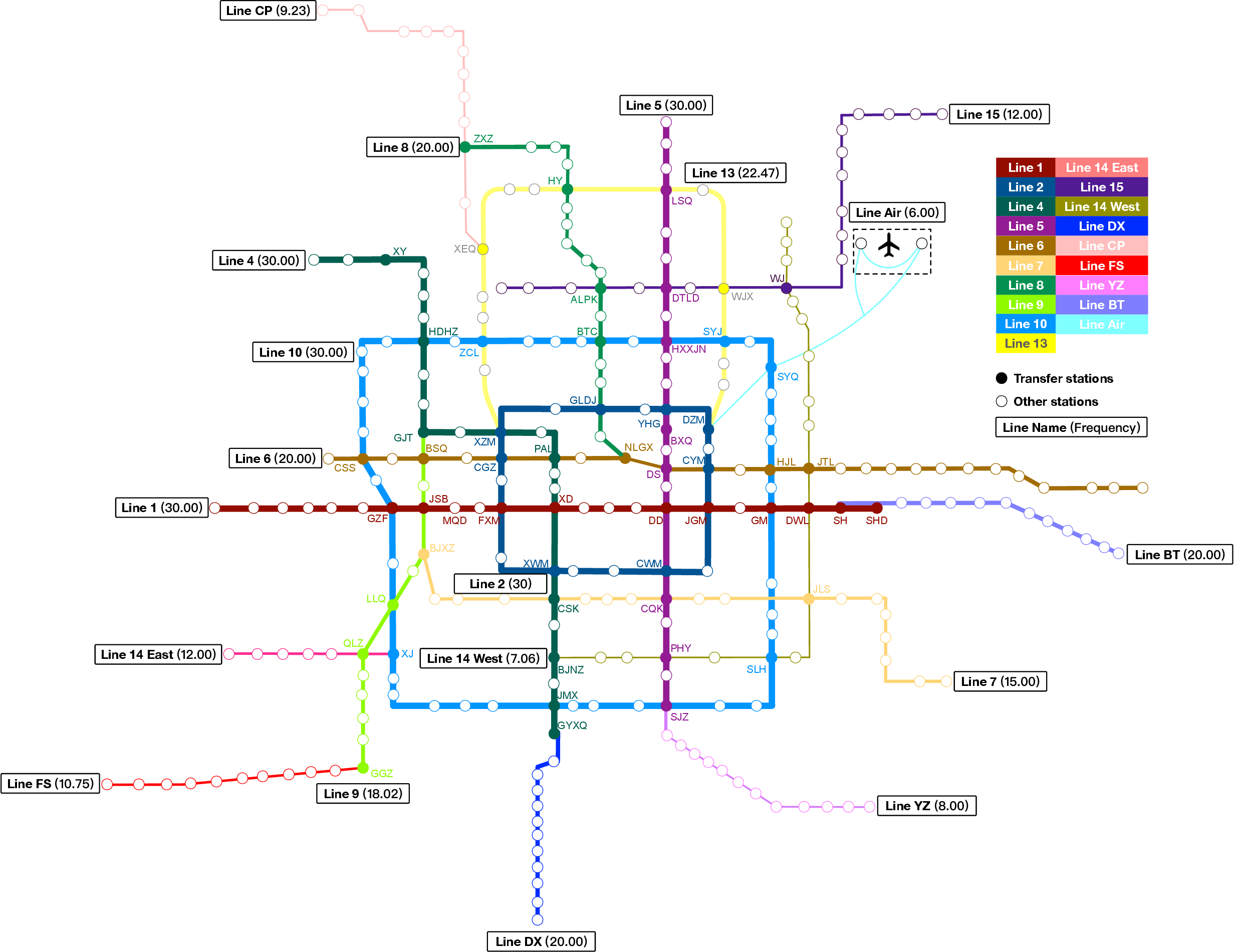}
         \caption{Existing line plan in the entire Beijing metro network}
     \end{subfigure}
    \caption{Real-life metro networks used in the computational study.}
    \label{fig_existingLineplan}
\end{figure}

The Beijing metro central sub-network includes lines 1, 2, 4, 5, 6, 8, 10, and 13 which cover the downtown area, as illustrated in Figure \ref{fig_existingLineplan}(a). For this sub-network, we focus on OD pairs where both the origin and destination lie within it. The Beijing metro sub-network and the entire network are of particular interest due to the large number of route choices available to commuters and the high utilization of the booming urban populations. In the Beijing metro network on April 2016, there are 277 stations (55 of which are transfer stations), and 19 existing lines, as shown in Figure \ref{fig_existingLineplan}(b). To build a line pool of the Beijing metro sub-network and the entire network, we first select the terminal stations of all the existing lines, and then select a subset of all stations including all the transfer stations, making sure that there is at least one route between any two of these terminal stations. After picking the stations use the shortest-path algorithm to find the paths with the shortest travel time between any two terminal stations, with each shortest path being a line. In this generating process, we still follow the restriction of established tracks, i.e., we can only set up lines if tracks already exist. For the line pool within the entire Beijing metro network, we further remove the dominated lines. The line pools obtained by the above method consists of 168 and 534 lines within the sub-network and the entire network, respectively. The transfer penalty is set as 15 minutes. The operational costs $d_\ell$ for line $\ell$ are equal to $N_\ell$ where $N_\ell$ is the number of stations traversed by line $\ell$. The set-up costs $e_{\ell}$ are again proportional to 
$d_\ell$, according to a factor that we vary in the experiments. We calculate the costs of the existing line plan used in real-life operations and set it as the baseline budget $\underline{B}$.

We select a working day to process the passenger demand, using OD data derived from the historically detected Automated Fare Collection system. On this day, the entire Beijing metro network accommodated a total of 5,989,805 passengers. During the busiest hour, between 7:00 a.m. and 8:00 a.m., 776,401 passengers traveled, experiencing severe crowding. The network comprises 56,916 OD pairs. Within the sub-network, there are 23,374 OD pairs and a passenger demand of 446,649. The running time on each section between two adjacent stations in this case study is set to the real value in real-world operations.

For the cut-and-column generation methods, we use the following parameters based on preliminary experiments. We use $\Delta_1=5$ and $\Delta_2=10$, i.e. if the average number of routes for each OD pair exceeds 5, we remove columns whose reduced costs exceed 10. In the first iteration of the column generation process, we terminate the process if the improvement within 20 consecutive iterations ($I = 20$) is not higher than 0.01\% ($\epsilon = 0.01 \%$). Subsequently, the values of $I$ and $\epsilon$ are dynamically adjusted: $I$ is decreased by 1, and $\epsilon$ is increased by 0.05\% in each successive iteration, with a minimum value of $I$ set to 1.

\subsection{Algorithmic performance}
\label{sec:GridResults}

We compare the performance of the proposed three procedures of the cut-and-column generation algorithm with GUROBI. To do so, we conduct two experiments using instances based on the grid network. In the first experiment, we use GUROBI to solve Formulation \eqref{formulation:SOconvexStrong} to optimality with varying numbers of routes for each OD pair, establishing a suitable benchmark for the algorithm. In the second set, we solve Formulation \eqref{formulation:SOconvexStrong} using three variants of the cut-and-column generation approach and compare the optimal solutions obtained from GUROBI.

\begin{figure}[htbp]
    \centering
\begin{tikzpicture}[scale=1]
\centering
\begin{axis}[xlabel={Number of routes in the PTN},
  ylabel={Passengers' perceived times},every axis plot/.append style={thick}, legend style={at={(0.05,0.95)},anchor=north west},legend cell align={left}]
\addplot[blue, mark=square*,mark options={fill=blue}] table [x=Number, y=Obj, col sep=comma] {gains.csv};
\end{axis}
\end{tikzpicture}
    \caption{Comparison of passengers' perceived travel times among various numbers of routes.}
    \label{fig_RoutesEffects}
\end{figure}
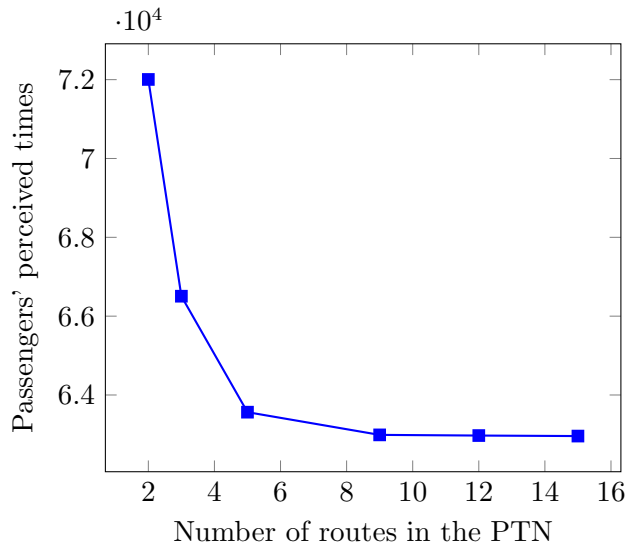
For the first experiment, we set a crowding factor of 0.05 and a budget of $\underline{B}$. We systematically increase the number of PTN routes from two to 15 for this set of experiments, solving all instances to optimality using GUROBI. The results are presented in Figure~\ref{fig:RoutesEffects}. We observe that if the number of routes is small, increasing the number of routes yields substantial reductions in passengers' perceived travel times. As the number of available routes increases, the marginal benefits taper off. These findings suggest that incorporating sufficient passengers' routes into line planning is necessary to model the routing adequately. Consequently, for subsequent experiments, we use GUROBI to solve the model with up to 15 PTN routes per OD pair, serving as a benchmark for the algorithm.

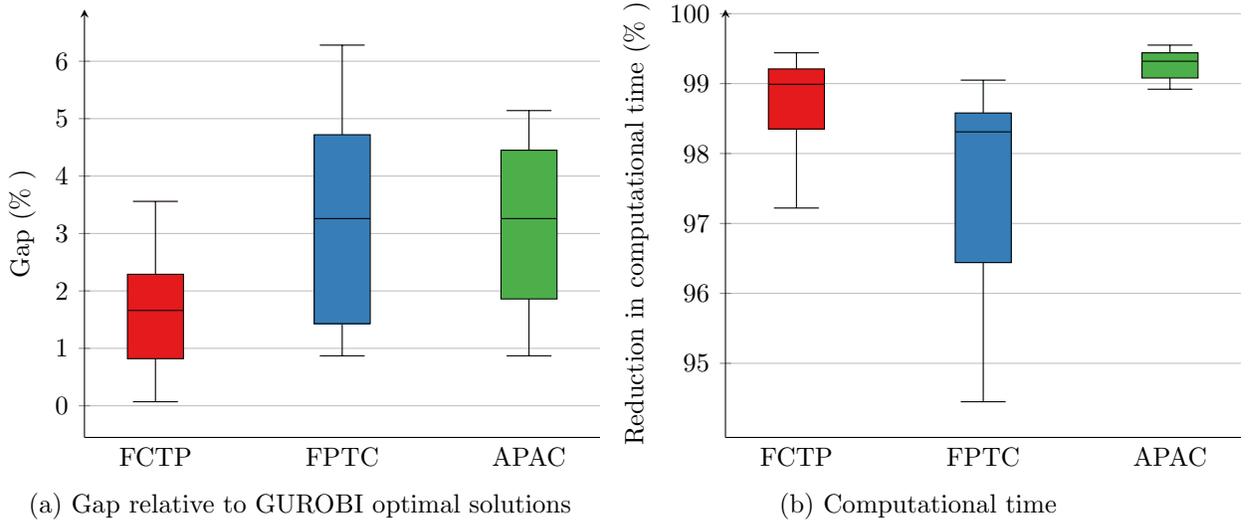
\begin{figure}[h]
    \centering
\begin{subfigure}[b]{0.49\textwidth}
    \centering
\begin{tikzpicture}
	\pgfplotstableread[col sep=comma]{Gap.csv}\csvdata
	\pgfplotstabletranspose\datatransposed{\csvdata} 
	\begin{axis}[
		boxplot/draw direction = y,
		x axis line style = {opacity=0},
		axis x line* = bottom,
		axis y line = left,
		enlarge y limits,
		ymajorgrids,
		xtick = {1, 2, 3},
		xticklabel style = {align=center, font=\small, rotate=0},
		xticklabels = {FCTP, FPTC, APAC},
        xticklabel style={font=\fontsize{10}{4}\selectfont}, 
		xtick style = {draw=none}, % Hide tick line
		ylabel = {Gap (\%
        )},
        ylabel style={font=\fontsize{10}{4}\selectfont},
		ytick = {0, 1, 2, 3, 4, 5, 6, 7, 8},
        yticklabel style={font=\fontsize{10}{4}\selectfont}, 
        boxplot/box extend=0.3
	]
		\foreach \n in {1,...,3} {
			\addplot+[boxplot, fill, draw=black] table[y index=\n] {\datatransposed};
		}
	\end{axis}
\end{tikzpicture}
\caption{Gap relative to GUROBI optimal solutions}
\end{subfigure}
\begin{subfigure}[b]{0.49\textwidth}
    \centering
\begin{tikzpicture}
	\pgfplotstableread[col sep=comma]{Time.csv}\csvdata
	\pgfplotstabletranspose\datatransposed{\csvdata} 
	\begin{axis}[
		boxplot/draw direction = y,
		x axis line style = {opacity=0},
		axis x line* = bottom,
		axis y line = left,
		enlarge y limits,
		ymajorgrids,
		xtick = {1, 2, 3},
		xticklabel style = {align=center, font=\small, rotate=0},
		xticklabels = {FCTP, FPTC, APAC},
        xticklabel style={font=\fontsize{10}{4}\selectfont}, 
		xtick style = {draw=none}, % Hide tick line
		ylabel = {Reduction in computational time (\%
        )},
        ylabel style={font=\fontsize{10}{4}\selectfont},
        yticklabel style={font=\fontsize{10}{4}\selectfont}, 
		ytick = {95, 96, 97, 98, 99, 100},
        boxplot/box extend=0.3
	]
		\foreach \n in {1,...,3} {
			\addplot+[boxplot, fill, draw=black] table[y index=\n] {\datatransposed};
		}
	\end{axis}
\end{tikzpicture}
\caption{Computational time}
\end{subfigure}
    \caption{Performance comparison of the three procedures versus GUROBI.}
    \label{fig_CGvsGUROBI}
\end{figure}

Next, Table~\ref{tab_algorithmVsGUROBI} and Figure~\ref{fig_CGvsGUROBI} provide a detailed performance comparison of the three variants of the cut-and-column generation algorithm (FCTP, FPTC, and APAC) against GUROBI. The table summarizes the objective values of the best feasible solutions produced by each algorithm, along with their computational times. It also includes the objective value achieved after 240 seconds of computation, the optimal objective value, and the corresponding computational time for GUROBI. The figure illustrates the average relative differences in objective values and computational time between the cut-and-column generation approach and GUROBI. Since we use the optimal solutions obtained by GUROBI as the benchmark, the relative gaps are true optimality gaps.

\begin{table}[h]
\centering
\caption{Comparison between three variants of the cut-and-column generation algorithm and GUROBI based on the $5 \times 5$-grid network.}
\label{tab_algorithmVsGUROBI}
\begin{threeparttable}
\resizebox{\textwidth}{!}{%
\begin{tabular}{ccccrcccccccccr}
\hline
\multicolumn{2}{c}{Instance}                                                                                      &  & \multicolumn{2}{c}{FCTP}                                        &  & \multicolumn{2}{c}{FPTC}                                        &  & \multicolumn{2}{c}{APAC}                                        &  & \multicolumn{3}{c}{GUROBI (15 PTN routes)} \\ \cline{1-2} \cline{4-5} \cline{7-8} \cline{10-11} \cline{13-15} 
\begin{tabular}[c]{@{}c@{}}Budget\\ (\%)\end{tabular} & \begin{tabular}[c]{@{}c@{}}Crowding\\ factor\end{tabular} &  & Objective & Time (s)&  & Objective & Time (s) &  & Objective & Time (s) &  & \begin{tabular}[c]{@{}c@{}}Objective\\ (240 s)\end{tabular} & \begin{tabular}[c]{@{}c@{}}Objective\\ (final)\end{tabular} & \begin{tabular}[c]{@{}c@{}}Time\\ (final, s)\end{tabular} \\
                                            \hline
100 & 0.01&  & 52,182 &  76  & & 52,862 & 130 & & 52,862 & 70 & & 52,969 & 52,147 & 13,677    \\
100 & 0.02&  & 55,295 & 76  & & 55,902 & 131 & & 55,902 & 60 && 55,520 & 55,116  &9,245  \\
100 & 0.03&  & 58,721 &  83  & & 58,942 & 150 & & 58,941 & 60 & &  59,201 & 57,866 & 8,873 \\
100 & 0.04&  & 61,780     & 97 & & 61,980 & 155 && 61,890 & 53 && 62,082 & 60,445 & 8,196   \\
100 & 0.05&  & 63,790 &  104  & & 65,012 & 152 && 65,012 & 66 & & 63,904  &  62,958 & 2,738  \\
100 & 0.06&  & 66,792   & 111 & & 68,046 & 186 & & 68,038 & 69 & & 66,288  & 65,465 & 8,601  \\
100 & 0.07&  &  69,524 &  132 & & 71,062 & 211 & & 71,062 & 64 & & 70,698 & 67,965 &  5,934   \\
100 & 0.08 & &  72,381  &  130 & & 74,084 & 210 && 74,082 & 81 & & 73,566 & 70,463  &  4,684   \\
\hline
120 & 0.01&  & 51,207  &  68 & & 52,351 & 133 && 52,351 & 85 & & 52,464 & 51,091  & 9,190 \\
120 & 0.02&  & 53,839 &  80 & & 54,820 & 129 && 54,820 & 63 & & 54,725 & 53,399  & 11,273   \\
120 & 0.03&  & 56,472 &   81 & & 57,289 & 127 && 57,289 & 55 & & 56,813 & 55,550  & 12,325    \\
120 & 0.04&  & 59,093&  93 & & 59,758 & 137& & 59,757 & 55 & & 59,902 & 57,602  & 11,474 \\
120 & 0.05&  & 60,951 &  96 & & 62,224 & 162 & & 62,224 & 55 & & 60,990 & 59,641  & 9,906  \\
120 & 0.06&  & 63,214 &  108& & 64,688 & 162 & & 64,688 & 68 & & 63,077 & 61,590  & 9,533   \\
120 & 0.07&  & 65,607 &  116 & & 67,150 & 200 & & 66,070 &  69 & & 65,161 & 63,540 & 14,120 \\
120 & 0.08&  & 66,950& 129  & & 69,405 & 194 && 68,409 & 70&  &  67,242 & 65,495  & 7,840  \\
\hline
\end{tabular}%
}
\begin{tablenotes}
    \footnotesize
    \item \textit{Note: Abbreviation: s = seconds.}
\end{tablenotes}
\end{threeparttable}
\end{table}

The results show that despite the relatively small size of these instances, GUROBI has a hard time finding optimal solutions, requiring multiple hours of computation time. The three variants of the cut-and-column generation approach find high-quality solutions in a matter of minutes, but with gaps up to about 5\%. FCTP achieves the smallest average optimality gap of 1.77\%, along with a narrower median and interquartile range compared to FPTC and APAC. Moreover, in 13 out of 16 instances, FCTP achieves better performance within at most 132 seconds of computation compared to GUROBI after 240 seconds. Based on these results and considering the strategic planning nature of our problem, we select FCTP as the solution method for further experimentation.

\subsection{Managerial insights}
\label{sec:managerialInsights}
This section evaluates the benefits of integrating crowding effects into the line planning decision and offers insights into optimized line plans across various crowding factors. The insights are derived by various experiments on the $5 \times 5$-grid network and the Beijing metro central sub-network. All models are solved using the FCTP variant of the cut-column generation algorithm. Our insights quantify the benefits of incorporating crowding effects into the line planning problem on reducing the crowding (Insight \ref{insightTraveltime}), on decreasing passengers' perceived travel times
(Insight~\ref{insightPerceived}), and on the structure of resulting line plans (Insight~\ref{insightLineplan}). In addition, our insights also quantify the benefits of integrating user-optimal routing in line planning (Insight~\ref{insightUE}).

\begin{insight}\label{insightTraveltime}
Incorporating crowding effects in line planning reduces crowding, even when the crowding factor is small. 
\end{insight}

We derive this insight through experiments on the Beijing metro central sub-network, focusing on two aspects. First, we optimize line plans by solving Formulation \eqref{formulation:SOconvexStrong} under varying crowding factors, considering line set-up costs of 30 and 40 with corresponding budget adjustments, to show the effectiveness of incorporating crowding effects in reducing crowding. Second, we compare the optimized line plans with the existing line plan used in real-life operations, highlighting the benefits of our approach in mitigating crowding compared to the practical benchmark.

\begin{figure}[h]
    \centering
% \begin{subfigure}[b]{0.49\textwidth}
%     \centering
% \begin{tikzpicture}[scale=1]
%     \centering
% \begin{axis}
% [xlabel={Congestion level (\%)},
%   ylabel={Cumulative distribution},every axis plot/.append style={thick}, legend style={font=\small, at={(0.33,0.30)},anchor=north west},legend cell align={left},  xmin=0, xmax=165 ]
% \addplot+[blue, dashed, mark=square, mark size=2pt, mark options={solid}] table [x=CongestionLevel, y=CF0001, col sep=comma] {gainsExisting.csv};
% \addlegendentry{Crowding factor = 0.001}
% \addplot+[Black, mark=o, mark size=2pt, mark options={solid}] table [x=CongestionLevel, y=CF0003, col sep=comma] {gainsExisting.csv};
% \addlegendentry{Crowding factor = 0.003}
% \addplot+[purple, mark=pentagon*, mark size=2pt, mark options={fill=purple}] table [x=CongestionLevel, y=CF0015, col sep=comma] {gainsExisting.csv};
% \addlegendentry{Crowding factor = 0.015}
% \end{axis}
% \end{tikzpicture}
%     \caption{Existing line plan}
% \end{subfigure}
\begin{subfigure}[b]{0.49\textwidth}
     \centering
\begin{tikzpicture}[scale=1]
    \centering
\begin{axis}
[xlabel={Crowding level (\%)},
  ylabel={Cumulative distribution},every axis plot/.append style={thick}, legend style={font=\small, at={(0.33,0.30)},anchor=north west},legend cell align={left},  xmin=0, xmax=165 ]
\addplot+[red, dashed, mark=square, mark size=2pt, mark options={solid}] table [x=CongestionLevel, y=CF0001, col sep=comma] {gainsRolf.csv};
\addlegendentry{Crowding factor = 0.001}
\addplot+[blue, mark=o, mark size=2pt, mark options={solid}] table [x=CongestionLevel, y=CF0003, col sep=comma] {gainsRolf.csv};
\addlegendentry{Crowding factor = 0.003}
\addplot+[YellowGreen, mark=triangle, mark size=3pt, mark options={fill=YellowGreen}, mark options={solid}] table [x=CongestionLevel, y=CF0015, col sep=comma] {gainsRolf.csv};
\addlegendentry{Crowding factor = 0.015}
\end{axis}
\end{tikzpicture}
\caption{Construction cost = 30}
\end{subfigure}
\begin{subfigure}[b]{0.49\textwidth}
     \centering
\begin{tikzpicture}[scale=1]
    \centering
\begin{axis}
[xlabel={Crowding level (\%)},
  ylabel={Cumulative distribution},every axis plot/.append style={thick}, legend style={font=\small, at={(0.33,0.30)},anchor=north west},legend cell align={left},  xmin=0, xmax=165 ]
\addplot+[red, dashed, mark=square, mark size=2pt, mark options={solid}] table [x=CongestionLevel, y=CF0001, col sep=comma] {gainsRolfE40.csv};
\addlegendentry{Crowding factor = 0.001}
\addplot+[blue, mark=o, mark size=2pt, mark options={solid}] table [x=CongestionLevel, y=CF0003, col sep=comma] {gainsRolfE40.csv};
\addlegendentry{Crowding factor = 0.003}
\addplot+[YellowGreen, mark=triangle, mark size=3pt, mark options={fill=YellowGreen}, mark options={solid}] table [x=CongestionLevel, y=CF0015, col sep=comma] {gainsRolfE40.csv};
\addlegendentry{Crowding factor = 0.015}
\end{axis}
\end{tikzpicture}
\caption{Construction cost = 40}
\end{subfigure}
\caption{Cumulative distribution of crowding levels under the optimized line plans with different crowding factors for the Beijing metro sub-network.}
\label{fig_CDFSubnetwork}
\end{figure}
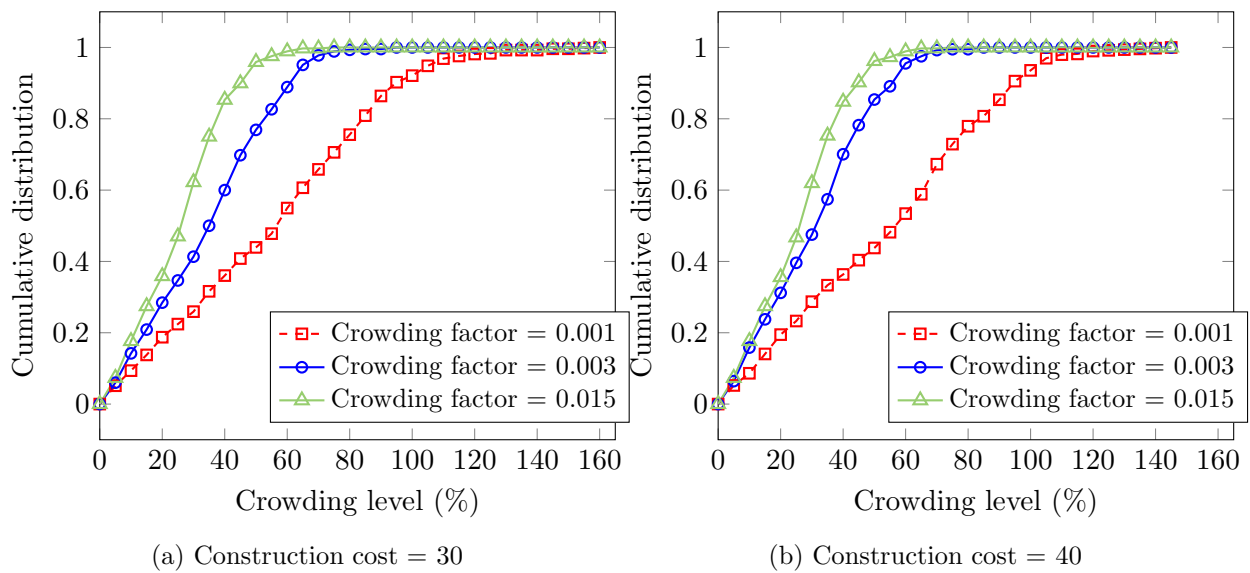

Figure~\ref{fig_CDFSubnetwork} presents the cumulative distribution of crowding levels under optimized line plans with various crowding
factors and different set-up costs. Regardless of the set-up costs, we can observe that as the crowding factor increases, the distribution shifts to the left, indicating that crowding is decreasing. For a crowding factor of 0.001, a significant portion of the arcs has a crowding level above 80\%. For the higher crowding factors, there are virtually no arcs with such high crowding levels. 
This suggests that a larger proportion of the network operates at lower crowding levels under higher crowding factors, improving overall passenger experience. Furthermore, comparing Figures~\ref{fig_CDFSubnetwork}(a) and (b), higher set-up costs not only reduce maximum crowing levels, but also compress the distribution towards lower crowding levels.

\begin{figure}[h]
    \centering
    \begin{tikzpicture}
    \begin{axis}[
        xlabel={Crowding factor},
        xtick={0.008,0.010, 0.012, 0.015},
        xticklabels={0.008, 0.010, 0.012, 0.015},
        ylabel={Relative reduction in APTT (\%)},
        ymin=4.5,
        ymax=8,
         width=8cm,
    height=8cm,
        enlargelimits=0.15,
        legend pos=north west, 
        ybar,
        bar width=7pt,
        scaled ticks=false,
    ]

        \addplot coordinates {
            (0.008,4.58) (0.010,4.64)
            (0.012,7.42) (0.015,7.43) 
        };
        \addplot coordinates {
            (0.008,7.37) (0.010,7.43)
            (0.012,7.44) (0.015,7.45) 
        };
        
        \legend{Installation costs $e_{\ell} = 30$, Installation costs $e_{\ell} = 40$}
    \end{axis}
    \end{tikzpicture}
\caption{Relative reductions in average perceived travel time (APTT) under optimized line plans with different crowding factors, compared to the existing real-life line plan, for instances on the Beijing metro sub-network.}
\label{fig_subnetworkResultsSOPercent}
\end{figure}
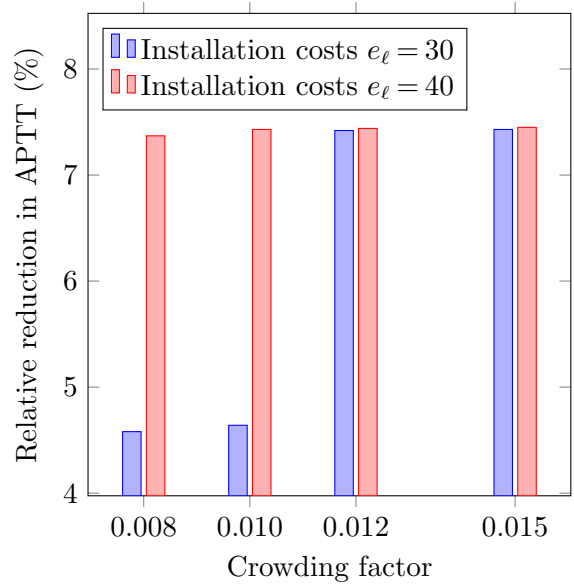

Next, Figure~\ref{fig_subnetworkResultsSOPercent} compares the existing line plan on the Beijing metros sub-network with optimized line plans. The relative reductions in average perceived travel time (APTT) is calculated as the relative differences between the results of the optimized and existing line plans. We observe that the optimized line plan reduces passengers' perceived travel time by 4.6\% to 7.5\%. Additionally, the results indicate that higher crowding factors correlate with greater reductions in the perceived travel time. Increasing set-up costs $e_\ell$ from 30 to 40 generally leads to further reductions in passengers' perceived travel time, although the differences become marginal at higher crowding factors. For instance, at a crowding factor of 0.010, the reductions of APTT are 4.6\% ($e_\ell=30$) and 7.4\% ($e_\ell=40$), respectively.

\begin{insight}\label{insightPerceived}
By incorporating crowding effects into line planning, passengers' perceived travel time is reduced. However, this comes at the expense of an increase in the travel time. 
\end{insight}

% 0.05, 2.56

\begin{figure}[h]
    \centering
\begin{subfigure}[b]{0.49\textwidth}
    \centering
    \begin{tikzpicture}
    \begin{axis}[
      xlabel={Crowding factor},
      xtick={0.008, 0.010, 0.012, 0.015},
    xticklabels={0.008, 0.010, 0.012, 0.015},
      ylabel={Relative reduction in APTT (\%)},
      width=8cm,
      height=8cm,
      xmin=0.008,
      ymin=0,
     scaled ticks=false,
      legend style={at={(0.5,-0.15)},anchor=north,legend columns=-1}
    ]
      \addplot[
        mark=square*,
        color=red
      ]
      table[col sep=comma, x=CongestionFactorsub, y=Improvesub] {gainsPTTSubnetwork.csv};
    \end{axis}
  \end{tikzpicture}
\caption{APTT }
\end{subfigure}
\begin{subfigure}[b]{0.49\textwidth}
    \centering
    \begin{tikzpicture}
    \begin{axis}[
      xlabel={Crowding factor},
      xtick={0.008, 0.010, 0.012, 0.015},
    xticklabels={0.008, 0.010, 0.012, 0.015},
      ylabel={Relative increase in ATT (\%)},
      width=8cm,
      height=8cm,
      xmin=0.008,
      ymin=0,
      scaled ticks=false,
      legend style={at={(0.5,-0.15)},anchor=north,legend columns=-1}
    ]
      \addplot+[blue, mark=*, mark size=2pt, mark options={fill=blue}, mark options={solid}]
      table[col sep=comma, x=CongestionFactorsub, y=Improvesub] {gainsTTSubnetwork.csv};
    \end{axis}
  \end{tikzpicture}
\caption{ATT}
\end{subfigure}
\caption{Relative differences in average perceived travel time (APTT) and average travel time (ATT) across various crowding factors compared to the benchmark without crowding effects under optimized line plans.}
\label{fig_PTTSO}
\end{figure}
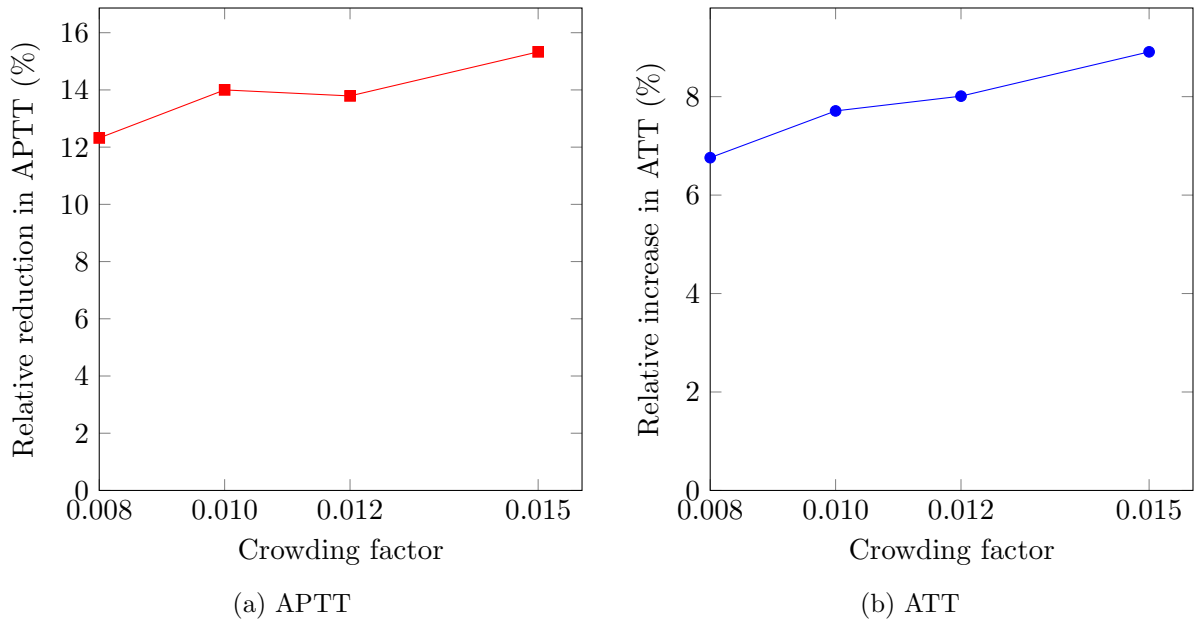

This insight is drawn from experiments on the Beijing metro central sub-network, using a benchmark that inputs the optimized line plan under crowding effects into Formulation~\eqref{formulation:NL} with $\gamma_a$ set to zero. For each benchmark, the perceived travel costs are recalculated by fixing the optimized line plan under crowding and applying the corresponding crowding factor values to $\gamma_a$. Thus, each point in Figure~\ref{fig_PTTSO} compares the perceived travel costs of passengers with and without crowding effects under the same line plan.

Figure~\ref{fig_PTTSO} presents the reduction in average perceived travel time (APTT) and the increase in average travel time (ATT). Three key observations emerge from these results. First, incorporating crowding effects reduces APTT, with reductions ranging from 12.3\% to 15.3\% under the presented crowding factor settings. This indicates that accounting for crowding effects in line planning enhances resource utilization, minimizing passengers' perceived travel time. Second, the incorporation of crowding effects increases ATT, as more passengers make transfers to search for routes with the smallest perceived travel costs, with increases ranging from 6.8\% to 8.9\%. Third, as the crowding factor increases, APTT decreases and ATT increases. However, the reduction in APTT is always greater than the increase in ATT, highlighting the overall benefit of considering crowding effects in line planning.

\begin{insight}\label{insightLineplan}Incorporating crowding effects results in structurally different line plans with fewer lines and higher frequencies. 
\end{insight}

\begin{figure}[h]
    \centering
\begin{subfigure}[b]{0.49\textwidth}
    \centering
   \includegraphics[height=0.6\textwidth]{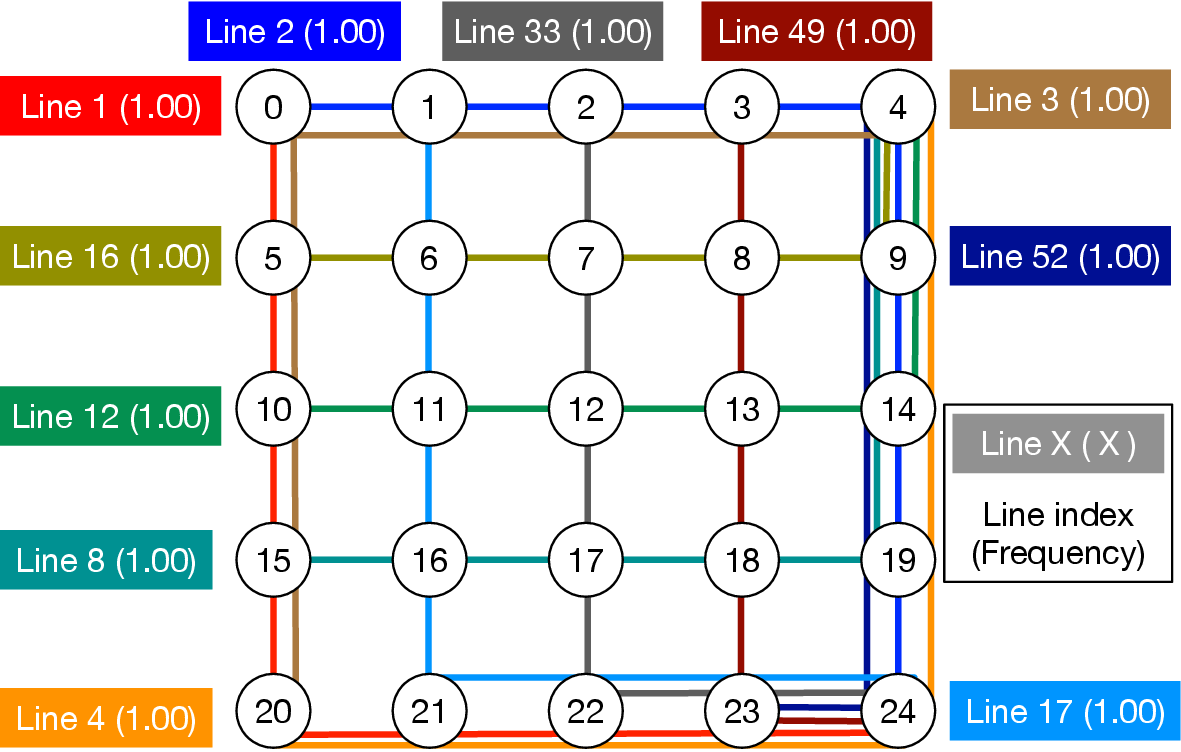}
    \caption{Crowding factor = 0}
    \end{subfigure}
\begin{subfigure}[b]{0.49\textwidth}
    \centering
\includegraphics[height=0.6\textwidth]{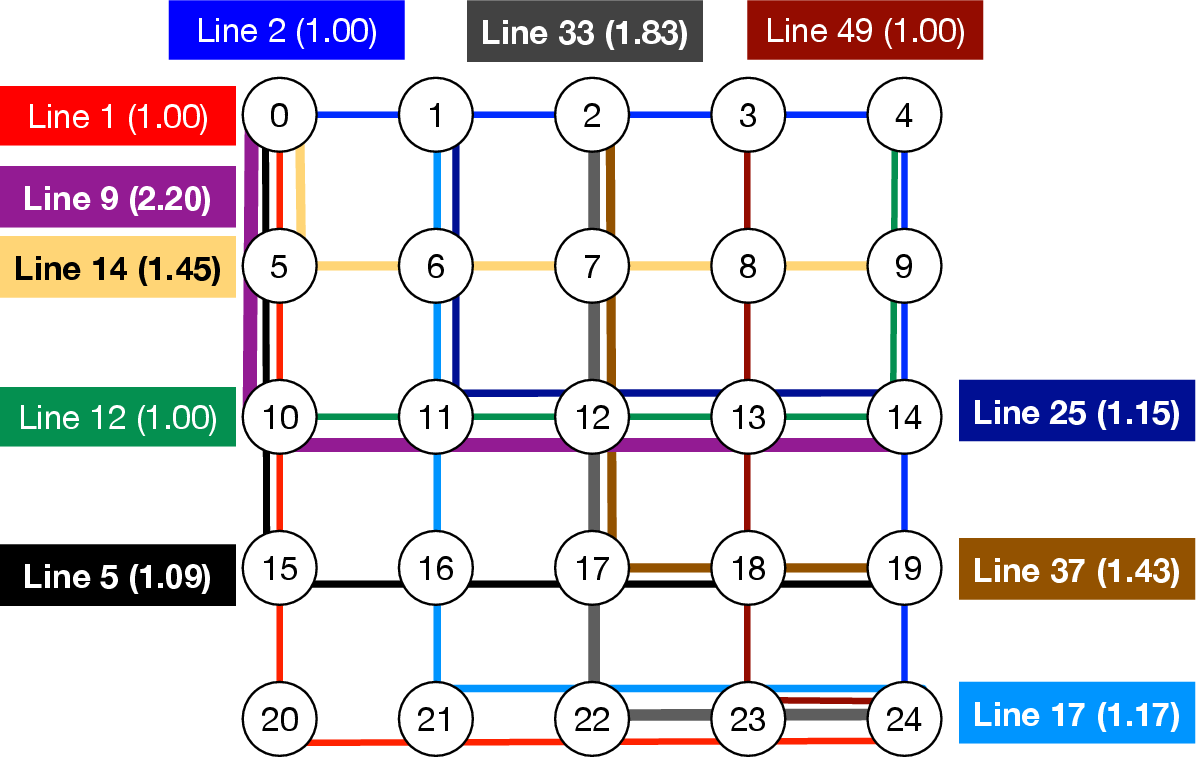}
    \caption{Crowding factor = 0.01}
    \end{subfigure}
    
\vspace{1em} 

\begin{subfigure}[b]{0.49\textwidth}
    \centering    
\includegraphics[height=0.6\textwidth]{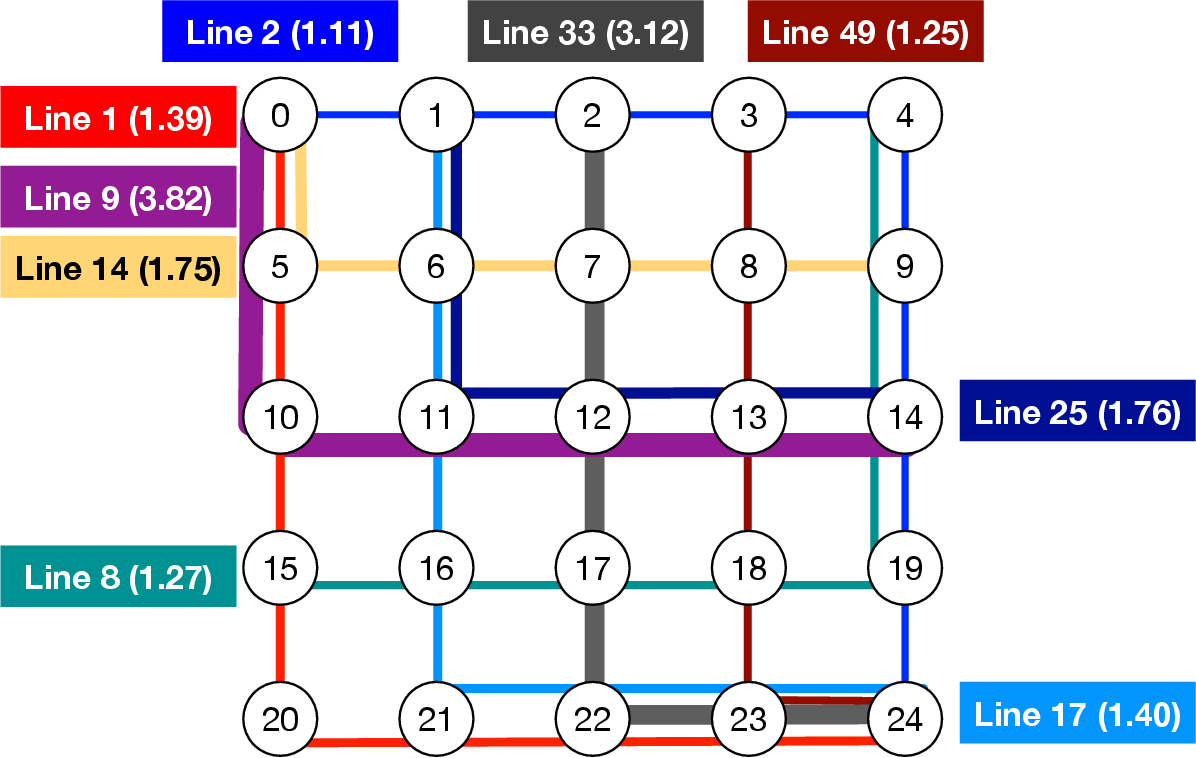}
    \caption{Crowding factor = 0.10}
    \end{subfigure}
\begin{subfigure}[b]{0.49\textwidth}
    \centering    
\includegraphics[height=0.6\textwidth]{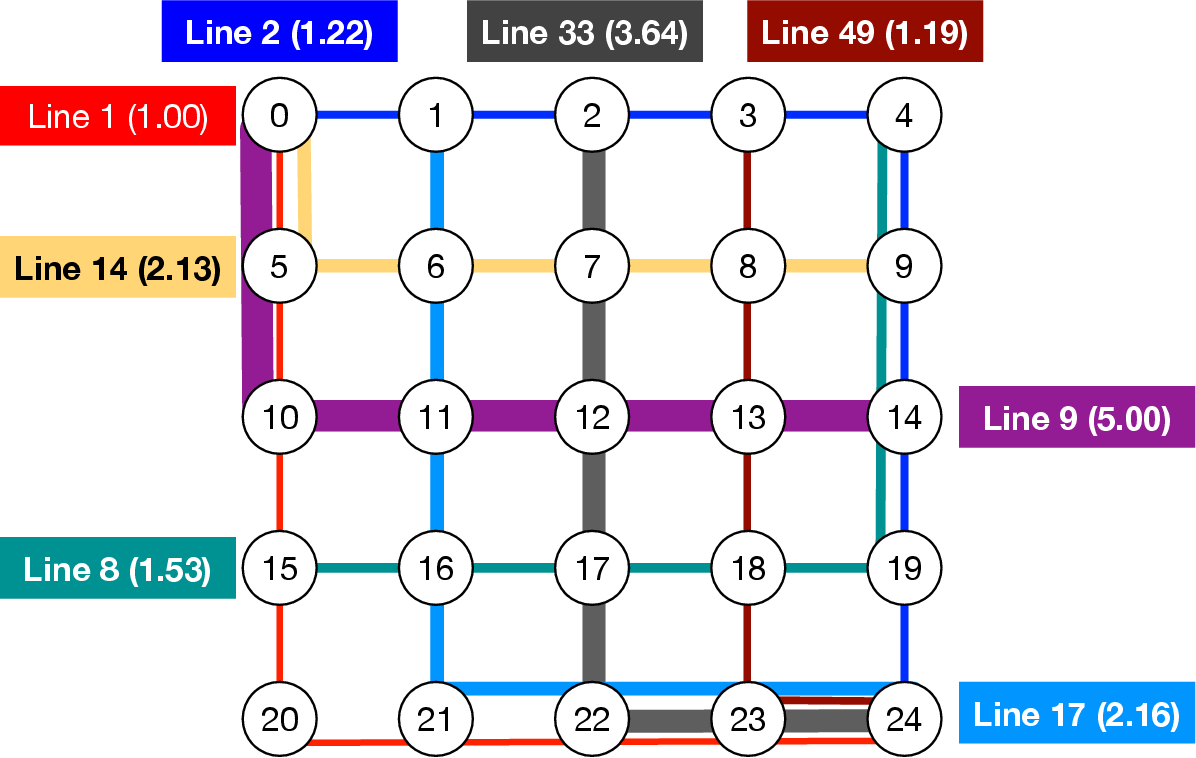}
    \caption{Crowding factor = 0.20}
    \end{subfigure}
    \caption{Optimized line plans for instances on the grid network with increasing values of the crowding factor.}
        \label{fig_B80LinePlan}
\end{figure}

In Figure \ref{fig_B80LinePlan}, we present the optimal line plans for the $5 \times 5$-grid network corresponding to the budget set at the baseline budget $\underline{B}$, highlighting both the operated lines and their frequencies. Without considering crowding effects, eleven lines are operated and all lines operate at the minimum frequency (i.e., 1). However, incorporating crowding penalties tends to reduce the number of lines while increasing their frequencies. For example, with a crowding factor of 0.01, eleven lines are operated, of which seven have frequencies greater than 1, and the maximum frequency used is 2.20. At a crowding factor of 0.20, only eight lines are built, only one of which has a frequency of 1, and the maximum frequency used increases to 5.

\begin{figure}[tbp]
    \centering
\begin{subfigure}[b]{\textwidth}
    \centering
   \includegraphics[height=0.6\textwidth]{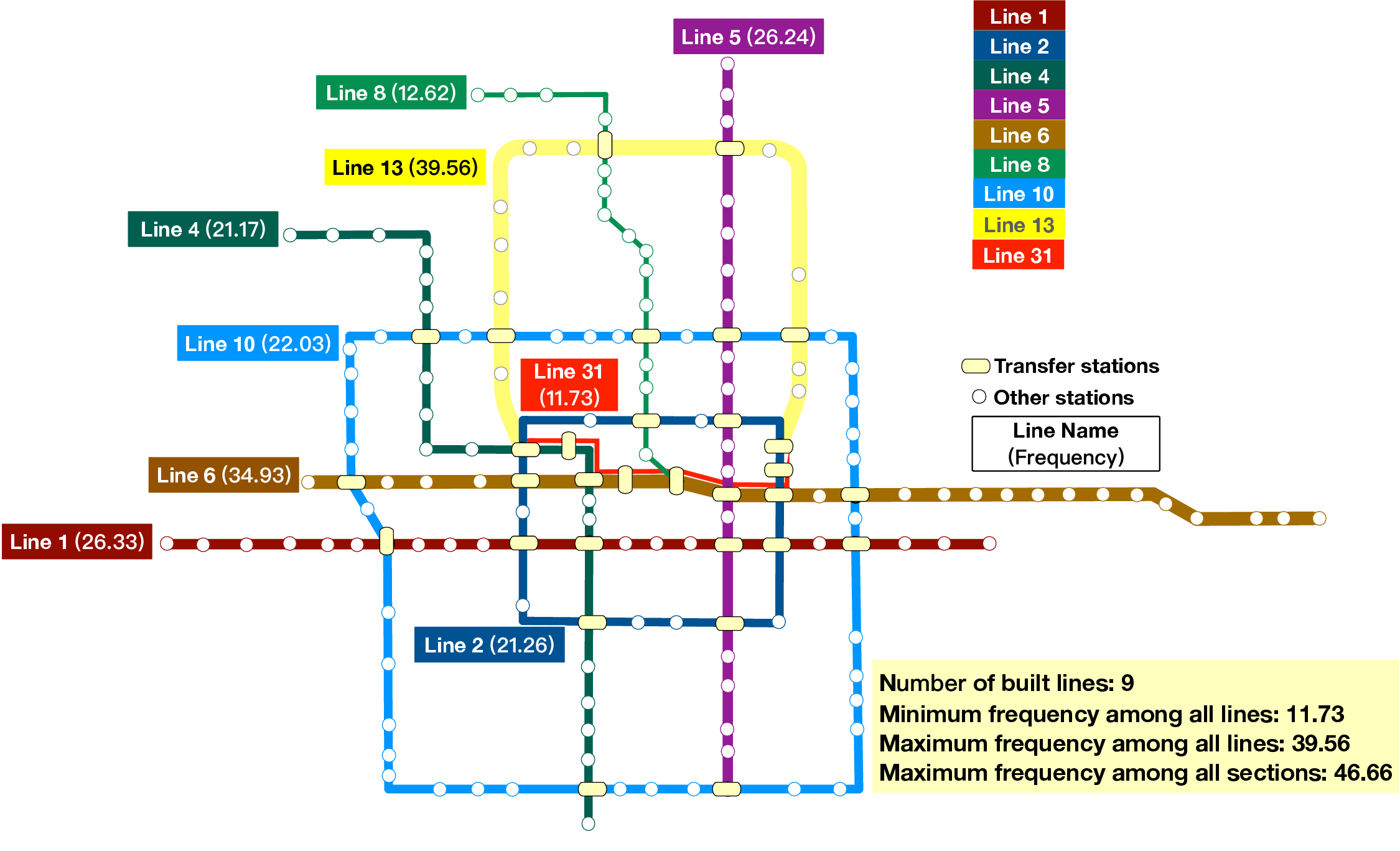}
    \caption{Crowding factor = 0.008}
    \end{subfigure}
    
\vspace{2em}

\begin{subfigure}[b]{\textwidth}
    \centering
\includegraphics[height=0.6\textwidth]{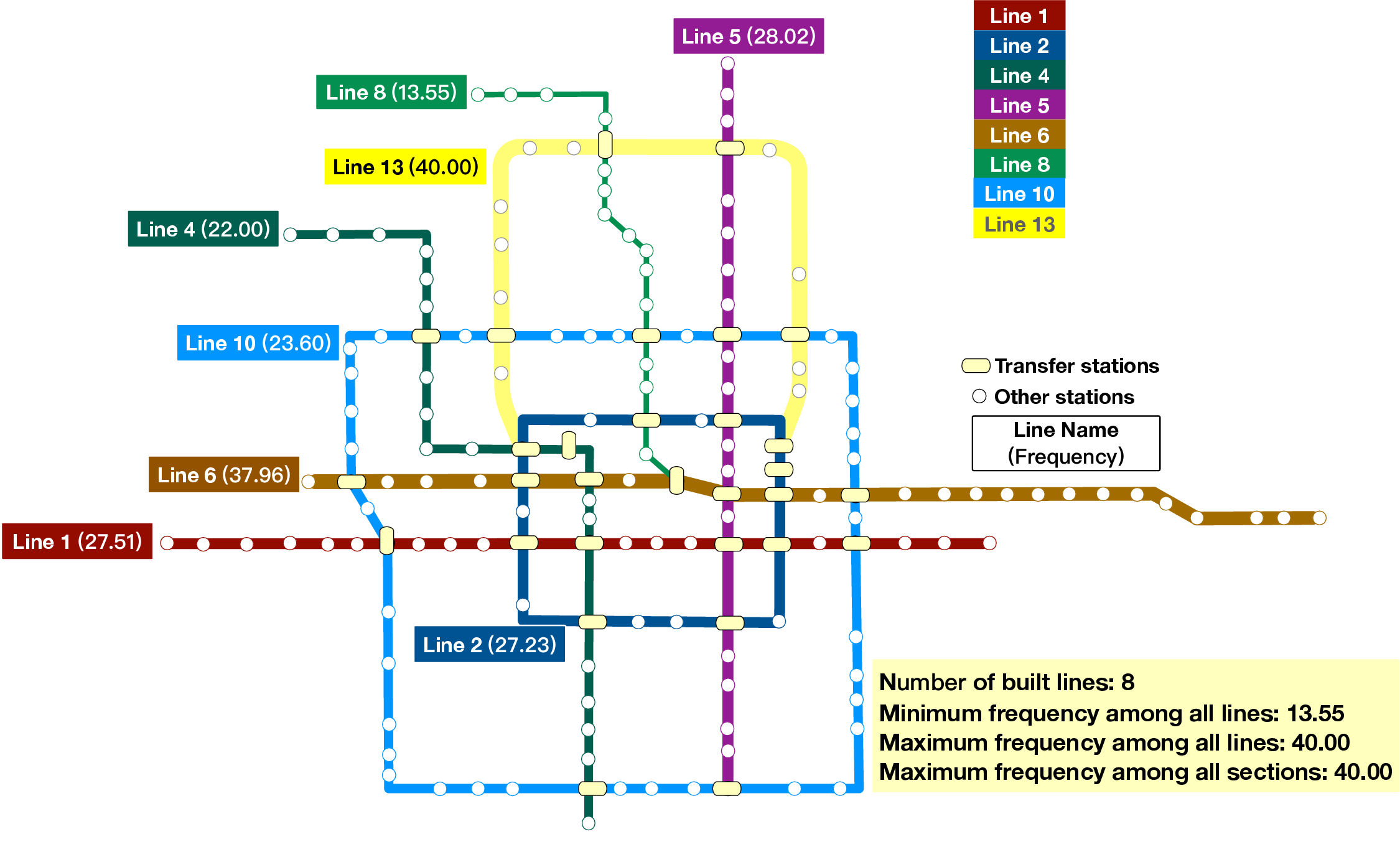}
    \caption{Crowding factor = 0.015}
    \end{subfigure}
    \caption{The optimal line plans for the Beijing metro sub-network with set-up costs set at 30.}
\label{fig_BeijingSubnetworkLinePlanE30}
\end{figure}

Furthermore, Figure \ref{fig_BeijingSubnetworkLinePlanE30} illustrates the optimized line plans for the Beijing metro sub-network under different crowding factors, with set-up costs set at 30. We observe that the optimized line plan contains 9 operated lines with a minimum headway frequency of 11.73 when the crowding factor is 0.008. As the crowding factor increases to 0.015, the line plan is reduced to 8 lines, with a minimum headway frequency of 13.55. A similar trend is observed for the optimized line plans with set-up costs set to 20, shown in Figure~\ref{fig_BeijingSubnetworkLineplanE20} in Appendix~\ref{sec:AppendixLinePlans}. These findings lead us to the conclusion that incorporating crowding effects in line planning results in fewer lines with higher frequencies. 

\begin{insight}\label{insightUE}
The system-optimal routing only deviates slightly from the user-equilibrium.
\end{insight}

Since passengers do not always resort to the recommended route through a network but
instead minimize their perceived travel time, we also compare it to the user equilibrium. To do so, we first determine an optimized line plan $\textbf{y}$ by solving Formulation~\eqref{formulation:SOconvexStrong}, with the budget set to $\underline{B}$, and re-evaluate passengers' perceived travel time using the cut-and-column generation algorithm. Next, we solve Formulation~\eqref{formulation:classicUE} using the same algorithm, utilizing the established line plan.

\begin{figure}[h]
    \centering
\begin{subfigure}[b]{0.49\textwidth}
\begin{tikzpicture}[scale=1]
\begin{axis}[xlabel={Crowding factor}, 
xtick={0.01, 0.02, 0.03, 0.04, 0.05, 0.06, 0.07, 0.08},
    xticklabels={0.01,0.02, 0.03, 0.04, 0.05, 0.06, 0.07, 0.08},
  ylabel={Relative difference between SO and UE (\%)},
  ytick={0, 0.01, 0.02, 0.03, 0.04, 0.05, 0.06, 0.07},
    yticklabels={0, 0.01, 0.02, 0.03, 0.04, 0.05, 0.06, 0.07},
  every axis plot/.append style={thick}, legend style={at={(0.05,0.95)},anchor=north west},legend cell align={left}, scaled ticks=false, xmin=0.01,     width=8cm, height=8cm]
\addplot [Black, mark=o, mark size=2pt, mark options={solid}] table [x=CongestionFactor, y=Dev, col sep=comma] {gainsUESOGrid.csv};
\end{axis}
\end{tikzpicture}
\caption{Grid network (budget $\underline{B})$}
\end{subfigure}
\begin{subfigure}[b]{0.49\textwidth}
    \centering
\begin{tikzpicture}[scale=1]
\centering
\begin{axis}[xlabel={Crowding factor},
xtick={0.008, 0.010, 0.012, 0.015},
    xticklabels={0.008, 0.010, 0.012, 0.015},
  ylabel={Relative difference between SO and UE (\%)},
   ytick={0.40, 0.45, 0.50, 0.55, 0.60, 0.65},
    yticklabels={0.40, 0.45, 0.50, 0.55, 0.60, 0.65},
    every axis plot/.append style={thick}, legend style={at={(0.05,0.95)},anchor=north west},legend cell align={left}, xmin=0.008, scaled ticks=false, 
 width=8cm, height=8cm]
\addplot [red, mark=square, mark size=2pt, mark options={solid}]table [x=CongestionFactor1, y=Dev1, col sep=comma] {gainsUESOBeijingSub.csv};
\end{axis}
\end{tikzpicture}
\caption{Beijing metro sub-network ($e_{\ell} = 40$)}
\end{subfigure}
\caption{Relative differences between UE and SO among various crowding penalties.}
    \label{fig_DiffUESO}
\end{figure}
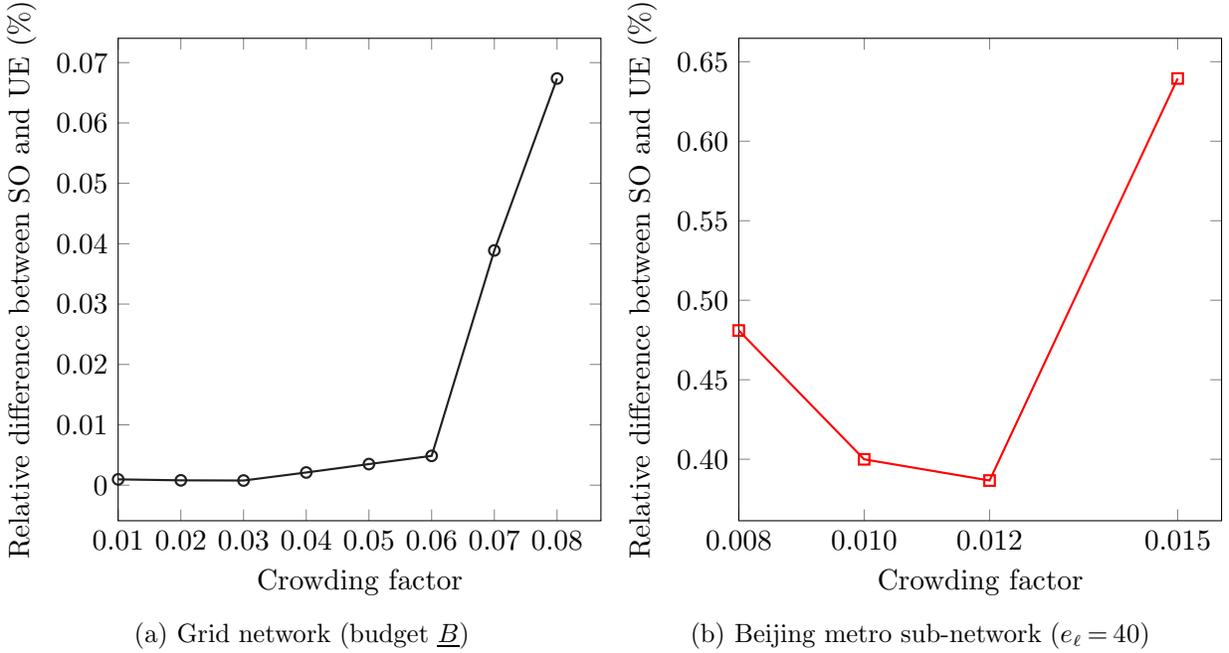

In Figure \ref{fig_DiffUESO}, we present the resulting relative differences between the UE and SO objective values, which are calculated by the optimal objective value of the UE model minus that of the SO model, and divide that of the SO model. From the results in Figure \ref{fig_DiffUESO}, a few observations stand out. First, we can observe that the user equilibria closely approximate the system optimum with respect to varying crowding factors in the instances based on the grid network. For example, the maximum value of the relative differences between the UE and SO objective values is below 0.07\% among all instances. A second observation is the divergence between UE and SO in the Beijing sub-network is a little more pronounced, with a maximum of 0.64\%. As
such, we conclude that integrating user-optimal routing in line planning offers limited benefits.

\subsection{Real-life case study based on the Beijing metro entire network}
\label{sec:entireNetwork}

We now assess the scalability of our algorithm through experiments on the entire Beijing metro network by varying the crowding factor from 0 to 0.008. The set-up costs $e_{\ell}$ are set to 40. Based on the conclusions in Section~\ref{sec:GridResults}, we use the FCTP variant of the cut-and-column generation algorithm. For the first column generation iteration, $\epsilon$ is started at 1\% and $I$ at 5. Then, with $I$ having a minimum value of 1, $\epsilon$ is raised by 0.2\% and $I$ is decreased by 1 in each iteration.

\begin{table}[h]
\centering
\caption{Computational times of the algorithm for instances on the entire Beijing metro network.}
\label{tab_entireTime}
\begin{threeparttable}
\resizebox{\textwidth}{!}{%
\begin{tabular}{rrrrr}
\hline
Crowding factor & Computational time (h) & Iterations (\#) & Average RMP time (s) & Average PP time (s) \\
\hline
0  &  2.2 & 53 &  10.8 & 132.5   \\
0.001  & 13.3 & 57 & 737.3  & 98.7  \\
0.003  & 33.7 & 71 & 1,597.4 & 106.0 \\
0.005  & 65.1 & 76 & 2,980.7 &  102.2\\         
0.008  & 108.4 & 78 & 4,895.3 & 107.1 \\
\hline
\end{tabular}%
}
\begin{tablenotes}
    \footnotesize
    \item \textit{Notes: Computational time includes both network generation and CPU time. Abbreviations: h = hours; s = seconds.}
\end{tablenotes}
\end{threeparttable}
\end{table}

Table~\ref{tab_entireTime} reports the results, including the total computational time, the number of iterations, the average time for solving the restricted master problem (RMP), and the average time for solving the pricing problem (PP). The instance with a crowding factor of 0, representing the line planning problem without crowding effects, is solved in 2.2 hours using our proposed cut-and-column generation algorithm. As the crowding factor increases, the computational time for solving the RMP rises significantly. For example, when the crowding factor increases from 0 to 0.001, the average computational time per RMP increases from 10.8 to 737.3 seconds. Similarly, when the crowding factor rises from 0.003 to 0.008, this time increases from 1,597.4 to 4,895.3 seconds. This growth is attributed to the increased computational challenge of solving the rotated second-order cone constraints \eqref{eq:tightLinedemand} by iteratively adding cutting planes. In contrast, the average computational time per PP does not exhibit a linear relationship with the incorporation of crowding effects. 

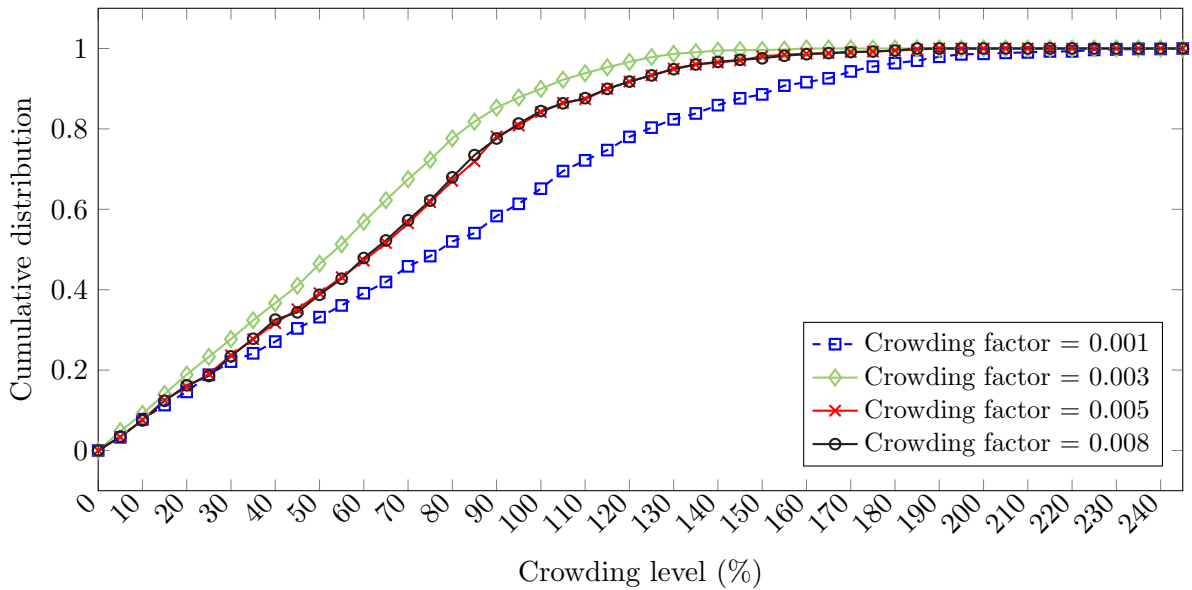
\begin{figure}
    \centering
\begin{tikzpicture}[scale=1]
    \centering
\begin{axis}
[xlabel={Crowding level (\%)},
  ylabel={Cumulative distribution},every axis plot/.append style={thick}, legend style={font=\small, at={(0.65,0.35)},anchor=north west},legend cell align={left},  xmin=0, xmax=245, xtick={0,10,...,240},
  xticklabel style={rotate=45, anchor=east},
    width=16cm,
    height=8cm ]
\addplot+[blue, dashed, mark=square, mark size=2pt, mark options={solid}] table [x=CongestionLevel, y=CF0001, col sep=comma] {gainsBeijingNetworkE40.csv};
\addlegendentry{Crowding factor = 0.001}
\addplot+[YellowGreen, mark=diamond, mark size=3pt, mark options={fill=YellowGreen}] table [x=CongestionLevel, y=CF0003, col sep=comma] {gainsBeijingNetworkE40.csv};
\addlegendentry{Crowding factor = 0.003}
\addplot+[red, mark=x, mark size=3pt, mark options={solid}] table [x=CongestionLevel, y=CF0005, col sep=comma] {gainsBeijingNetworkE40.csv};
\addlegendentry{Crowding factor = 0.005}
\addplot+[Black, mark=o, mark size=2pt, mark options={solid}] table [x=CongestionLevel, y=CF0008, col sep=comma] {gainsBeijingNetworkE40.csv};
\addlegendentry{Crowding factor = 0.008}
\end{axis}
\end{tikzpicture}
\caption{Cumulative distribution of crowding levels under the optimized line plans under different crowding factors for the Beijing metro entire network.}
\label{fig_CDFBeijingNetwork}
\end{figure}

Figure \ref{fig_CDFBeijingNetwork} shows the cumulative distribution of crowding levels under different crowding factors for the Beijing metro entire network. We can observe that as the crowding factor increases from 0.001 to 0.003, the cumulative distribution shifts to the left. This observation aligns with Insight \ref{insightTraveltime}, which suggests that incorporating crowding effects into line planning reduces crowding. This reduction enhances passenger comfort and improves the efficiency of the metro system. These findings highlight the importance of incorporating crowding effects in line planning for large-scale networks like the Beijing metro. However, for larger crowding factors (0.005 and 0.008), the problem becomes increasingly complex, and the diving heuristic has a harder time finding good solutions, resulting in solutions with more crowding. This shows that practitioners should carefully select the crowding factor.

\section{Conclusion}
\label{sec:conclusion}

This paper addressed a line planning problem under crowding, extending traditional line planning by explicitly incorporating crowding effects into passengers' routing decisions. The problem is formulated as a mixed-integer nonlinear programming model. To obtain its computationally tractable form, we propose an MI-SOCP reformulation with second-order cone constraints. To solve the MI-SOCP, we develop a cut-and-column generation methodology. Specifically, a cut generation procedure is designed to dynamically generate linear cutting planes to address the rotated second-order cone constraints on-the-fly, while a column generation procedure is implemented to iteratively identify effective passenger routing. Within this framework, we propose three variants of the cut-and-column generation approach.

Extensive computational experiments yield multiple takeaways. First, the cut-and-column generation algorithm delivers high-quality solutions with an average optimality gap ranging from 1.8\% to 3.5\% on the $5 \times 5$ grid network within reasonable computation times. The proposed methodology scales effectively to large-scale, real-world problems involving over 200 stations, 530 candidate lines, and 56,000 OD pairs, producing optimized line plans that outperform the one currently used in practical operations. Additionally,
our proposed methodology provides strong practical benefits, achieving significant reductions in passengers' perceived travel times while reducing crowding within the network. At a time when increasingly crowded public transit systems need high-quality services, these benefits offer solid support and valuable tools for practitioners to optimize resource planning effectively.

\section*{Acknowledgments}
This work was supported by the National Natural Science Foundation of China (Grant No. 72288101).

\bibliographystyle{informs2014trsc} 
\bibliography{sample}

\begin{thebibliography}{39}
\providecommand{\natexlab}[1]{#1}
\providecommand{\url}[1]{\texttt{#1}}
\providecommand{\urlprefix}{URL }

\bibitem[{Basciftci \protect\BIBand{} Van~Hentenryck(2023)}]{Basciftci2023}
Basciftci B, Van~Hentenryck P, 2023 \emph{Capturing travel mode adoption in designing on-demand multimodal transit systems}. \emph{Transportation Science} 57(2):351--375.

\bibitem[{Bayram, Y\i{}ld\i{}z, \protect\BIBand{} Farham(2023)}]{Vedat2023}
Bayram V, Y\i{}ld\i{}z B, Farham MS, 2023 \emph{Hub network design problem with capacity, congestion, and stochastic demand considerations}. \emph{Transportation Science} 57(5):1276--1295.

\bibitem[{{Beijing Municipal Government Press Office}(2023)}]{beijing2023subway}
{Beijing Municipal Government Press Office}, 2023 \emph{New technology breaks bottlenecks: Frequency on the beijing metro is expected to decrease to 1 minute and 30 seconds.} \urlprefix\url{https://www.peopleapp.com/rmharticle/30002598532}, accessed: 2024-07-01.

\bibitem[{Bertsimas, Ng, \protect\BIBand{} Yan(2021)}]{Bertsimas2021}
Bertsimas D, Ng YS, Yan J, 2021 \emph{Data-driven transit network design at scale}. \emph{Operations Research} 69(4):1118--1133.

\bibitem[{Bornd\"{o}rfer, Gr\"{o}tschel, \protect\BIBand{} Pfetsch(2007)}]{Borndorfer2007}
Bornd\"{o}rfer R, Gr\"{o}tschel M, Pfetsch ME, 2007 \emph{A column-generation approach to line planning in public transport}. \emph{Transportation Science} 41(1):123--132.

\bibitem[{Cacchiani, Qi, \protect\BIBand{} Yang(2020)}]{Cacchiani2020}
Cacchiani V, Qi J, Yang L, 2020 \emph{Robust optimization models for integrated train stop planning and timetabling with passenger demand uncertainty}. \emph{Transportation Research Part B: Methodological} 136:1--29.

\bibitem[{Canca et~al.(2019)Canca, De-Los-Santos, Laporte, \protect\BIBand{} Mesa}]{Canca2019}
Canca D, De-Los-Santos A, Laporte G, Mesa JA, 2019 \emph{Integrated railway rapid transit network design and line planning problem with maximum profit}. \emph{Transportation Research Part E: Logistics and Transportation Review} 127:1--30.

\bibitem[{Chai et~al.(2024)Chai, Yin, D’Ariano, Liu, Yang, \protect\BIBand{} Tang}]{Chai2024}
Chai S, Yin J, D’Ariano A, Liu R, Yang L, Tang T, 2024 \emph{A branch-and-cut algorithm for scheduling train platoons in urban rail networks}. \emph{Transportation Research Part B: Methodological} 181:102891.

\bibitem[{Daganzo(1997)}]{Daganzo1997}
Daganzo CF, 1997 \emph{Fundamentals of Transportation and Traffic Operations} (Emerald Publishing).

\bibitem[{{De Vos}, {Van Lieshout}, \protect\BIBand{} Dollevoet(2024)}]{Marelot2024}
{De Vos} MH, {Van Lieshout} RN, Dollevoet T, 2024 \emph{Electric vehicle scheduling in public transit with capacitated charging stations}. \emph{Transportation Science} 58(2):279--294.

\bibitem[{Franceschetti et~al.(2018)Franceschetti, Honhon, Laporte, \protect\BIBand{} Van~Woensel}]{Franceschetti2018}
Franceschetti A, Honhon D, Laporte G, Van~Woensel T, 2018 \emph{A shortest-path algorithm for the departure time and speed optimization problem}. \emph{Transportation Science} 52(4):756--768.

\bibitem[{Friedrich et~al.(2017{\natexlab{a}})Friedrich, Hartl, Schiewe, \protect\BIBand{} Sch{\"o}bel}]{GridRef}
Friedrich M, Hartl M, Schiewe A, Sch{\"o}bel A, 2017{\natexlab{a}} \emph{Angebotsplanung im {\"o}ffentlichen verkehr-planerische und algorithmische l{\"o}sungen}. \emph{Heureka’17} 9.

\bibitem[{Friedrich et~al.(2017{\natexlab{b}})Friedrich, Hartl, Schiewe, \protect\BIBand{} Sch\"{o}bel}]{SchobelAnita2017}
Friedrich M, Hartl M, Schiewe A, Sch\"{o}bel A, 2017{\natexlab{b}} \emph{{Integrating passengers' assignment in cost-optimal line planning}}. \emph{17th Workshop on Algorithmic Approaches for Transportation Modelling, Optimization, and Systems (ATMOS 2017)}, volume~59 of \emph{Open Access Series in Informatics (OASIcs)}, 5:1--5:16 (Dagstuhl, Germany: Schloss Dagstuhl -- Leibniz-Zentrum f{\"u}r Informatik).

\bibitem[{Gattermann, Harbering, \protect\BIBand{} Sch{\"o}bel(2017)}]{Gattermann2017}
Gattermann P, Harbering J, Sch{\"o}bel A, 2017 \emph{Line pool generation}. \emph{Public Transport} 9:7--32.

\bibitem[{Goerigk \protect\BIBand{} Schmidt(2017)}]{Goerigk2017}
Goerigk M, Schmidt M, 2017 \emph{Line planning with user-optimal route choice}. \emph{European Journal of Operational Research} 259(2):424--436.

\bibitem[{Goossens, {Van Hoesel}, \protect\BIBand{} Kroon(2004)}]{Goossens2004}
Goossens JW, {Van Hoesel} S, Kroon L, 2004 \emph{A branch-and-cut approach for solving railway line-planning problems}. \emph{Transportation Science} 38(3):379--393.

\bibitem[{Goossens, {Van Hoesel}, \protect\BIBand{} Kroon(2006)}]{Goossens2006}
Goossens JW, {Van Hoesel} S, Kroon L, 2006 \emph{On solving multi-type railway line planning problems}. \emph{European Journal of Operational Research} 168(2):403--424.

\bibitem[{Hartleb et~al.(2023)Hartleb, Schmidt, Huisman, \protect\BIBand{} Friedrich}]{Hartleb2023}
Hartleb J, Schmidt M, Huisman D, Friedrich M, 2023 \emph{Modeling and solving line planning with mode choice}. \emph{Transportation Science} 57(2):336--350.

\bibitem[{Jiang, Rasmussen, \protect\BIBand{} Nielsen(2022)}]{Jiang2022}
Jiang Y, Rasmussen TK, Nielsen OA, 2022 \emph{Integrated optimization of transit networks with schedule- and frequency-based services subject to the bounded stochastic user equilibrium}. \emph{Transportation Science} 56(6):1452--1468.

\bibitem[{K\i{}nay, Gzara, \protect\BIBand{} Alumur(2023)}]{Burak2023}
K\i{}nay OB, Gzara F, Alumur SA, 2023 \emph{Charging station location and sizing for electric vehicles under congestion}. \emph{Transportation Science} 57(6):1433--1451.

\bibitem[{Kowalczyk \protect\BIBand{} Leus(2018)}]{Kowalczyk2018}
Kowalczyk D, Leus R, 2018 \emph{A branch-and-price algorithm for parallel machine scheduling using zdds and generic branching}. \emph{INFORMS Journal on Computing} 30(4):768--782.

\bibitem[{Lu et~al.(2023)Lu, Yang, Yang, Zhou, \protect\BIBand{} Gao}]{Lu2023}
Lu Y, Yang L, Yang H, Zhou H, Gao Z, 2023 \emph{Robust collaborative passenger flow control on a congested metro line: A joint optimization with train timetabling}. \emph{Transportation Research Part B: Methodological} 168:27--55.

\bibitem[{Lu et~al.(2022)Lu, Yang, Yang, Gao, Zhou, Meng, \protect\BIBand{} Qi}]{Lu2022}
Lu Y, Yang L, Yang K, Gao Z, Zhou H, Meng F, Qi J, 2022 \emph{A distributionally robust optimization method for passenger flow control strategy and train scheduling on an urban rail transit line}. \emph{Engineering} 12:202--220.

\bibitem[{Luan \protect\BIBand{} Corman(2022)}]{Luan2022}
Luan X, Corman F, 2022 \emph{Passenger-oriented traffic control for rail networks: An optimization model considering crowding effects on passenger choices and train operations}. \emph{Transportation Research Part B: Methodological} 158:239--272.

\bibitem[{Roughgarden(2007)}]{roughgarden2007}
Roughgarden T, 2007 \emph{Routing games}. Nisan N, Roughgarden T, Tardos E, Vazirani VV, eds., \emph{Algorithmic game theory}, chapter~18, 459--484 (Cambridge University Press).

\bibitem[{Schiewe, Schiewe, \protect\BIBand{} Schmidt(2019)}]{Alexander2019}
Schiewe A, Schiewe P, Schmidt M, 2019 \emph{The line planning routing game}. \emph{European Journal of Operational Research} 274(2):560--573.

\bibitem[{Schiewe, J\"ager, \protect\BIBand{} Sch\"{o}bel(2024)}]{schiewe2024lintim}
Schiewe P, J\"ager S, Sch\"{o}bel A, 2024 \emph{Lintim - integrated optimization in public transportation}. \url{https://www.lintim.net/}, open source.

\bibitem[{Schiewe \protect\BIBand{} Sch\"{o}bel(2020)}]{Schiewe2020}
Schiewe P, Sch\"{o}bel A, 2020 \emph{Periodic timetabling with integrated routing: Toward applicable approaches}. \emph{Transportation Science} 54(6):1714--1731.

\bibitem[{Schmidt(2012)}]{Schmidt2012}
Schmidt M, 2012 \emph{Line planning with equilibrium routing}. Ph.D. thesis, Institute for Numerical and Applied Mathematics, University of Gottingen.

\bibitem[{Schmidt \protect\BIBand{} Schöbel(2024)}]{schmidt2024planning}
Schmidt M, Schöbel A, 2024 \emph{Planning and optimizing transit lines}. \urlprefix\url{https://arxiv.org/abs/2405.10074}.

\bibitem[{Sch\"{o}bel(2012)}]{Schobel2012}
Sch\"{o}bel A, 2012 \emph{Line planning in public transportation: models and methods}. \emph{OR Spectrum} 34:491--510.

\bibitem[{Sch{\"o}bel(2017)}]{Schobel2017}
Sch{\"o}bel A, 2017 \emph{An eigenmodel for iterative line planning, timetabling and vehicle scheduling in public transportation}. \emph{Transportation Research Part C: Emerging Technologies} 74:348--365.

\bibitem[{Sch{\"o}bel \protect\BIBand{} Scholl(2005)}]{Schobel2005}
Sch{\"o}bel A, Scholl S, 2005 \emph{Line planning with minimal transfers}. \emph{Proceeding of 5th Workshop on Algorithmic Methods and Models for Optimization of Railways} (Germany).

\bibitem[{{United Nations}(2019)}]{un75}
{United Nations}, 2019 \emph{{Shifting Demographics: A Visual Guide}}. \url{https://www.un.org/en/un75/shifting-demographics}, [Accessed March 28, 2023].

\bibitem[{{Van Lieshout}, Bouman, \protect\BIBand{} Huisman(2020)}]{Lieshout2020}
{Van Lieshout} RN, Bouman PC, Huisman D, 2020 \emph{Determining and evaluating alternative line plans in out-of-control situations}. \emph{Transportation Science} 54(3):740--761.

\bibitem[{Wardrop(1952)}]{Wardrop1952}
Wardrop JG, 1952 \emph{Some theoretical aspects of road traffic research}. \emph{Proceedings of the Institute of Civil Engineers} Part II(1):325--378.

\bibitem[{Xia, Ma, \protect\BIBand{} {Sharif Azadeh}(2024)}]{Xia2024integrated}
Xia D, Ma J, {Sharif Azadeh} S, 2024 \emph{Integrated timetabling, vehicle scheduling, and dynamic capacity allocation of modular autonomous vehicles under demand uncertainty}. \urlprefix\url{https://arxiv.org/abs/2410.16409}.

\bibitem[{Yao, Nie, \protect\BIBand{} Fu(2024)}]{Yao2024}
Yao Z, Nie L, Fu H, 2024 \emph{Railway line planning with passenger routing: Direct-service network representations and a two-phase solution approach}. \emph{Transportation Research Part B: Methodological} 186:102989.

\bibitem[{Yin et~al.(2021)Yin, D’Ariano, Wang, Yang, \protect\BIBand{} Tang}]{Yin2021}
Yin J, D’Ariano A, Wang Y, Yang L, Tang T, 2021 \emph{Timetable coordination in a rail transit network with time-dependent passenger demand}. \emph{European Journal of Operational Research} 295(1):183--202.

\end{thebibliography}

\newpage 

\begin{APPENDICES}

\section{Optimized line plans for the Beijing metro central sub-network} \label{sec:AppendixLinePlans}

This section presents the optimized line plans for the Beijing metro central sub-network, with set-up costs set to 20, as shown in Figure~\ref{fig_BeijingSubnetworkLineplanE20}.

\begin{figure}[h]
    \centering
\begin{subfigure}[b]{\textwidth}
    \centering
   \includegraphics[height=0.55\textwidth]{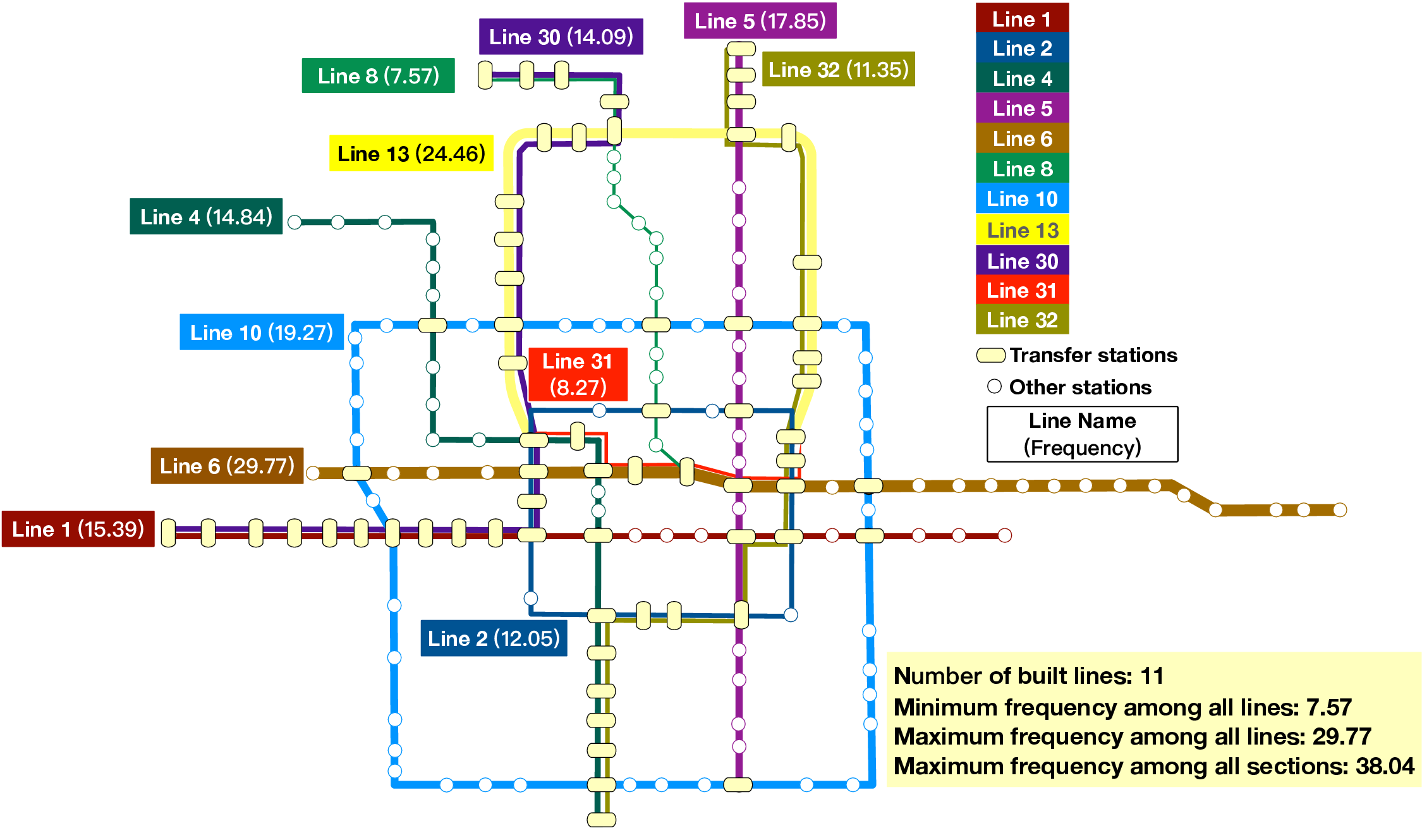}
    \caption{Crowding factor = 0.008}
    \end{subfigure}
    
\vspace{1em}

\begin{subfigure}[b]{\textwidth}
    \centering
\includegraphics[height=0.55\textwidth]{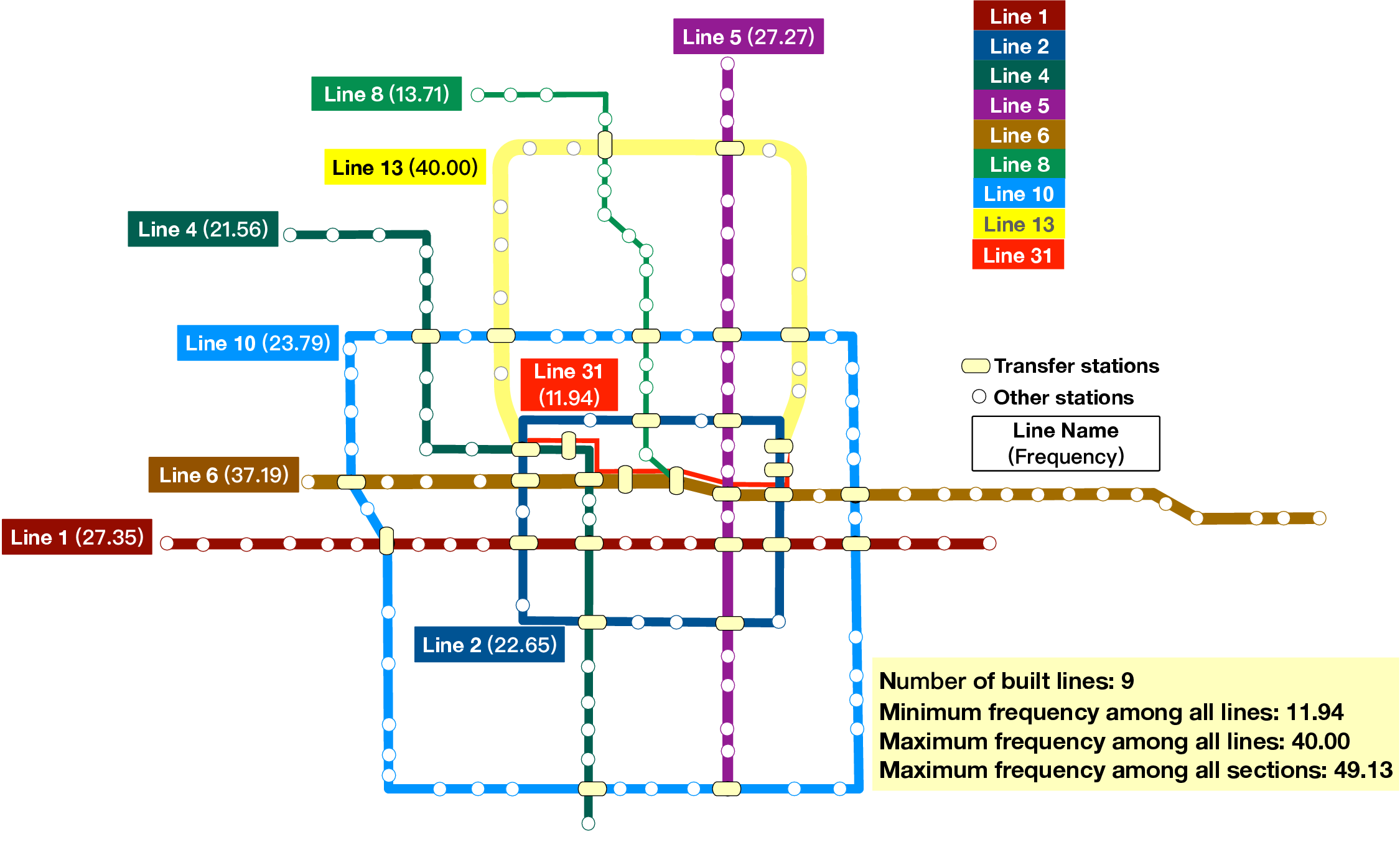}
    \caption{Crowding factor = 0.015}
    \end{subfigure}
    \caption{The optimal line plans for the Beijing metro sub-network with set-up costs set at 20.}
\label{fig_BeijingSubnetworkLineplanE20}
\end{figure}

 \end{APPENDICES}

\end{document}